\documentclass[11pt]{amsart}
\usepackage[T1]{fontenc}
\usepackage[margin=3cm]{geometry}
\usepackage{amsmath,amssymb,amsthm,mathtools}
\usepackage[hidelinks]{hyperref}
\usepackage{tikz}
\numberwithin{equation}{section}
\newtheorem{theorem}{Theorem}[section]
\newtheorem{proposition}[theorem]{Proposition}
\newtheorem{lemma}[theorem]{Lemma}
\theoremstyle{remark}
\newtheorem*{remark}{Remark}
\newcommand{\R}{\mathbb R}
\newcommand{\dd}{\,\mathrm d}
\title[Prescribing finite first singular sets]{Prescribing finite first singular sets of mean curvature flow in \texorpdfstring{$\R^3$}{R^3}}
\author{Jingwen Chen}
\address{Department of Mathematics, Faculty of Science,
University of Macau,\newline
E11, Avenida da Universidade, Taipa,
Macau}
\email{jingwenchen@um.edu.mo}
\date{}
\begin{document}
\begin{abstract}
For any $n$ distinct points in $\R^3$ and any sufficiently small positive time, we construct a smooth mean curvature flow of embedded mean-convex tori whose first singularities occur precisely at those points at that time.
\end{abstract}
\maketitle

\section{Introduction}\label{sec:introduction}
A smooth family of surfaces $M_t\subset\R^3$ moves by mean curvature flow if
\[
 (\partial_t x)^\perp=\vec H,
\]
where $\vec H$ is the mean curvature vector. A mean curvature flow starting from a smooth closed surface develops singularities in finite time.

White \cite{White2000} proved that the spacetime singular set of mean curvature flow from a closed mean-convex surface in $\R^3$ has parabolic Hausdorff dimension at most one. The perturbation results of Chodosh--Choi--Mantoulidis--Schulze \cite{CCMS2024} and Chodosh--Choi--Schulze \cite{CCS2023}, together with the multiplicity-one theorem of Bamler--Kleiner \cite{BK2023}, imply that a generic closed embedded surface in $\R^3$ evolves with only spherical and cylindrical singularities. For mean curvature flows of closed embedded surfaces with only spherical and cylindrical singularities, Colding--Minicozzi \cite{ColdingMinicozzi2016} show that the cylindrical singular set is contained in finitely many compact embedded Lipschitz curves. Ao Sun, Zhihan Wang, and Jinxin Xue \cite{SWXRegularity2025} prove that the degenerate cylindrical singular set is locally contained in $C^{2,\alpha}$ curves in spacetime for every $0<\alpha<1$. Daniels-Holgate--Hershkovits--Sz{\'e}kelyhidi \cite{DanielsHolgateHershkovitsSzekelyhidi2026} prove that mean curvature flow from a generic compact initial surface in $\R^3$ has only finitely many spacetime singular points.

The corresponding realization problem asks which sets can be prescribed as the entire first singular set of a mean curvature flow from a connected initial surface. The basic examples in $\R^3$ are a single point for the shrinking sphere and a line for the shrinking cylinder. In the marriage ring example discussed in \cite{White2002}, a thin torus of revolution collapses to a circle under mean curvature flow. Yiqi Huang and Xinrui Zhao \cite{HuangZhao2025} construct non-product incomplete mean curvature flows whose first singular set is a line. For closed connected initial surfaces, few constructions prescribe several first singular points. For each $g\geq3$, Ilmanen--White \cite{IlmanenWhite2025} used $g$-fold rotational symmetry to construct a closed embedded surface of genus $g$ whose first singular set consists of the vertices of a regular $g$-gon.

Ao Sun and Jinxin Xue \cite{SXII} record a conjecture of Colding--Minicozzi that any smooth curve in $\R^3$ can be the singular set of a mean-convex mean curvature flow. Haslhofer \cite[Problem~5.6]{HaslhoferICM2026} also asks which closed curves can occur as singular sets of mean curvature flows of embedded surfaces and notes that the only known example is the marriage ring. Sun--Xue \cite{SXII} showed that if the first singular set of a smooth embedded mean curvature flow is a closed $C^{2,\alpha}$ curve of cylindrical singularities, an arbitrarily small $C^{2,\alpha}$ normal perturbation of the initial surface can produce a single first singular point near any chosen point of the curve.

Tsiamis \cite{Tsiamis2026} constructs a noncompact mean-convex ancient mean curvature flow in $\R^{d+m}$, $m\geq2$, whose first singular set is $K\times\{0\}$ for any nonempty closed set $K\subset\R^d$, using a smooth ambient metric arbitrarily $C^\infty$-close to the Euclidean metric. We work in Euclidean $\R^3$ with smooth closed connected initial surfaces.

\begin{theorem}\label{thm:main}
For distinct points $p_1,\ldots,p_n\in\R^3$, there is $T_0=T_0(p_1,\ldots,p_n)>0$ such that, for every $0<T<T_0$, there is a smooth embedded strictly mean-convex torus $M_0\subset\R^3$ whose mean curvature flow first becomes singular at time $T$, precisely at $p_1,\ldots,p_n$.
\end{theorem}

\begin{remark}
Each prescribed singularity is a nondegenerate neckpinch in the sense of \cite{SXI,SXII}.
\end{remark}

To the best of the author's knowledge, this is the first result prescribing an arbitrary nonempty finite subset of $\R^3$ as the first singular set of a mean curvature flow starting from a smooth closed connected surface. Here the initial surface is a torus, independently of the number or configuration of the prescribed points.

The problem of prescribing the singular set has also been studied for other geometric flows. For harmonic map heat flow into $\mathbb S^2$, D\'avila, del Pino, and Juncheng Wei \cite{DavilaDelPinoWei2020} realized any nonempty finite set of interior points as the entire blow-up set in a smooth bounded planar domain. Together with Pesce \cite{DavilaDelPinoPesceWei2019}, they constructed solutions on smooth bounded axially symmetric domains in $\R^3$ that blow up precisely on a prescribed circle about the symmetry axis. For K\"ahler--Ricci flow, Jian Song and Weinkove \cite{SongWeinkove2013} proved that the fibers over finitely many blown-up points of a compact K\"ahler manifold collapse simultaneously to their respective points at the first singular time if the limiting cohomology class is the pullback of a K\"ahler class on the original manifold.

Sun--Xue \cite{SXII} showed that a nondegenerate neckpinch at the first singular time of a flow from a closed embedded hypersurface persists nearby in spacetime under sufficiently small $C^2$ perturbations of the initial hypersurface. We track how the singular points and times vary with the initial torus to make all its necks pinch at the prescribed points at time $T$. Sun--Wang--Xue \cite{SWXPassing2025} proved that each nondegenerate neckpinch is isolated in spacetime.

\subsection*{Outline of the proof}
\begin{figure}[htbp]
\centering
\begin{tikzpicture}[x=.70cm,y=.70cm,line cap=round,line join=round,font=\small]
\def\constructioncurve{
 (-2.2187,-0.5797) (-2.3124,-0.6031) (-2.4062,-0.6266) (-2.5000,-0.6500) (-2.5938,-0.6734) (-2.6876,-0.6969) (-2.7813,-0.7203) (-2.8351,-0.7330)
 (-2.8869,-0.7407) (-2.9349,-0.7396) (-2.9779,-0.7270) (-3.0146,-0.7013) (-3.0444,-0.6615) (-3.0666,-0.6076) (-3.0810,-0.5404) (-3.0874,-0.4612)
 (-3.0860,-0.3716) (-3.0769,-0.2739) (-3.0606,-0.1704) (-3.0375,-0.0638) (-3.0084,0.0431) (-2.9737,0.1479) (-2.9343,0.2479) (-2.8909,0.3409)
 (-2.8442,0.4249) (-2.7949,0.4984) (-2.7438,0.5604) (-2.6914,0.6107) (-2.6382,0.6496) (-2.5848,0.6786) (-2.5312,0.6998) (-2.4778,0.7167)
 (-2.3852,0.7444) (-2.2926,0.7722) (-2.2000,0.8000) (-2.1074,0.8278) (-2.0148,0.8556) (-1.9222,0.8833) (-1.7968,0.9214) (-1.6710,0.9617)
 (-1.5444,1.0061) (-1.4169,1.0555) (-1.2882,1.1106) (-1.1584,1.1712) (-1.0273,1.2370) (-0.8952,1.3072) (-0.7619,1.3806) (-0.6278,1.4559)
 (-0.4930,1.5316) (-0.3577,1.6060) (-0.2221,1.6775) (-0.0864,1.7443) (0.0491,1.8050) (0.1842,1.8580) (0.3188,1.9021) (0.4527,1.9363)
 (0.5857,1.9601) (0.7178,1.9732) (0.8489,1.9758) (0.9791,1.9688) (1.1085,1.9535) (1.2373,1.9319) (1.3656,1.9069) (1.4604,1.8879)
 (1.5552,1.8690) (1.6500,1.8500) (1.7448,1.8310) (1.8396,1.8121) (1.9344,1.7931) (2.0471,1.7691) (2.1599,1.7367) (2.2723,1.6896)
 (2.3837,1.6233) (2.4933,1.5353) (2.5999,1.4244) (2.7023,1.2910) (2.7990,1.1366) (2.8885,0.9637) (2.9693,0.7758) (3.0399,0.5768)
 (3.0989,0.3713) (3.1449,0.1641) (3.1767,-0.0401) (3.1933,-0.2364) (3.1940,-0.4204) (3.1783,-0.5882) (3.1461,-0.7363) (3.0975,-0.8623)
 (3.0334,-0.9649) (2.9549,-1.0437) (2.8636,-1.1000) (2.7617,-1.1366) (2.6522,-1.1582) (2.5386,-1.1711) (2.4424,-1.1808) (2.3462,-1.1904)
 (2.2500,-1.2000) (2.1538,-1.2096) (2.0576,-1.2192) (1.9614,-1.2289) (1.8727,-1.2370) (1.7847,-1.2413) (1.6981,-1.2386) (1.6134,-1.2267)
 (1.5308,-1.2041) (1.4505,-1.1703) (1.3725,-1.1250) (1.2966,-1.0687) (1.2227,-1.0024) (1.1505,-0.9273) (1.0795,-0.8451) (1.0095,-0.7575)
 (0.9400,-0.6667) (0.8705,-0.5745) (0.8006,-0.4832) (0.7300,-0.3945) (0.6582,-0.3104) (0.5849,-0.2323) (0.5100,-0.1615) (0.4331,-0.0988)
 (0.3543,-0.0444) (0.2737,0.0018) (0.1915,0.0408) (0.1079,0.0741) (0.0237,0.1042) (-0.0675,0.1361) (-0.1588,0.1681) (-0.2500,0.2000)
 (-0.3412,0.2319) (-0.4325,0.2639) (-0.5237,0.2958) (-0.5921,0.3191) (-0.6604,0.3389) (-0.7285,0.3523) (-0.7963,0.3573) (-0.8639,0.3527)
 (-0.9312,0.3378) (-0.9981,0.3124) (-1.0649,0.2768) (-1.1314,0.2320) (-1.1979,0.1790) (-1.2642,0.1191) (-1.3305,0.0542) (-1.3970,-0.0142)
 (-1.4635,-0.0841) (-1.5303,-0.1538) (-1.5974,-0.2215) (-1.6649,-0.2856) (-1.7327,-0.3447) (-1.8009,-0.3977) (-1.8696,-0.4439) (-1.9388,-0.4829)
 (-2.0083,-0.5149) (-2.0782,-0.5407) (-2.1484,-0.5616)}
\def\constructionouter{
 (-2.1949,-0.6747) (-2.2919,-0.6854) (-2.3874,-0.7016) (-2.4818,-0.7228) (-2.5750,-0.7485) (-2.6670,-0.7792) (-2.7576,-0.8154) (-2.8134,-0.8393)
 (-2.8777,-0.8612) (-2.9538,-0.8729) (-3.0420,-0.8621) (-3.1313,-0.8202) (-3.2082,-0.7525) (-3.2693,-0.6689) (-3.3077,-0.5726) (-3.3173,-0.4679)
 (-3.3156,-0.3586) (-3.3051,-0.2450) (-3.2867,-0.1280) (-3.2610,-0.0092) (-3.2286,0.1095) (-3.1900,0.2261) (-3.1457,0.3386) (-3.0958,0.4452)
 (-3.0406,0.5445) (-2.9751,0.6319) (-2.8870,0.6927) (-2.7973,0.7365) (-2.7133,0.7685) (-2.6360,0.7898) (-2.5666,0.8026) (-2.5060,0.8105)
 (-2.4096,0.8257) (-2.3148,0.8464) (-2.2216,0.8718) (-2.1297,0.9019) (-2.0392,0.9368) (-1.9504,0.9772) (-1.8345,1.0429) (-1.7260,1.1264)
 (-1.6223,1.2174) (-1.5037,1.2685) (-1.3822,1.3205) (-1.2588,1.3781) (-1.1331,1.4413) (-1.0048,1.5094) (-0.8740,1.5814) (-0.7407,1.6563)
 (-0.6050,1.7325) (-0.4670,1.8084) (-0.3269,1.8822) (-0.1846,1.9524) (-0.0402,2.0169) (0.1061,2.0743) (0.2542,2.1228) (0.4037,2.1610)
 (0.5541,2.1879) (0.7042,2.2028) (0.8531,2.2058) (0.9988,2.1918) (1.1334,2.1246) (1.2604,2.0570) (1.3848,2.0029) (1.4771,1.9711)
 (1.5704,1.9449) (1.6647,1.9235) (1.7600,1.9069) (1.8562,1.8953) (1.9536,1.8892) (2.0758,1.8886) (2.2125,1.8892) (2.3672,1.8792)
 (2.5149,1.8123) (2.6488,1.7048) (2.7747,1.5740) (2.8915,1.4218) (2.9989,1.2503) (3.0965,1.0619) (3.1835,0.8596) (3.2590,0.6470)
 (3.3218,0.4280) (3.3709,0.2068) (3.4051,-0.0124) (3.4231,-0.2259) (3.4238,-0.4302) (3.4057,-0.6226) (3.3667,-0.8012) (3.3038,-0.9641)
 (3.2139,-1.1074) (3.0973,-1.2242) (2.9508,-1.2844) (2.8025,-1.2902) (2.6696,-1.2796) (2.5484,-1.2686) (2.4508,-1.2652) (2.3539,-1.2674)
 (2.2575,-1.2746) (2.1615,-1.2866) (2.0661,-1.3036) (1.9712,-1.3263) (1.8815,-1.3530) (1.7866,-1.3827) (1.6839,-1.4120) (1.5714,-1.4325)
 (1.4558,-1.4213) (1.3476,-1.3760) (1.2456,-1.3168) (1.1507,-1.2465) (1.0625,-1.1674) (0.9801,-1.0818) (0.9023,-0.9916) (0.8280,-0.8989)
 (0.7565,-0.8054) (0.6870,-0.7132) (0.6191,-0.6243) (0.5523,-0.5405) (0.4865,-0.4635) (0.4218,-0.3944) (0.3579,-0.3341) (0.2948,-0.2826)
 (0.2396,-0.2266) (0.1922,-0.1549) (0.1351,-0.0897) (0.0673,-0.0349) (-0.0087,0.0117) (-0.0955,0.0561) (-0.1843,0.0950) (-0.2748,0.1292)
 (-0.3668,0.1589) (-0.4605,0.1838) (-0.5561,0.2033) (-0.6266,0.2123) (-0.6920,0.2119) (-0.7498,0.1999) (-0.7974,0.1763) (-0.8340,0.1452)
 (-0.8653,0.1190) (-0.9025,0.1032) (-0.9458,0.0801) (-0.9947,0.0471) (-1.0485,0.0040) (-1.1063,-0.0481) (-1.1672,-0.1078) (-1.2308,-0.1733)
 (-1.2969,-0.2427) (-1.3653,-0.3140) (-1.4362,-0.3855) (-1.5097,-0.4554) (-1.5862,-0.5221) (-1.6659,-0.5839) (-1.7542,-0.6307) (-1.8532,-0.6503)
 (-1.9486,-0.6598) (-2.0371,-0.6656) (-2.1189,-0.6700)}
\def\constructioninner{
 (-2.2424,-0.4846) (-2.3330,-0.5208) (-2.4250,-0.5515) (-2.5182,-0.5772) (-2.6126,-0.5984) (-2.7081,-0.6146) (-2.8051,-0.6253) (-2.8569,-0.6267)
 (-2.8961,-0.6201) (-2.9161,-0.6064) (-2.9138,-0.5920) (-2.8980,-0.5824) (-2.8806,-0.5704) (-2.8639,-0.5464) (-2.8542,-0.5083) (-2.8575,-0.4545)
 (-2.8563,-0.3846) (-2.8487,-0.3027) (-2.8345,-0.2128) (-2.8141,-0.1185) (-2.7881,-0.0232) (-2.7574,0.0697) (-2.7229,0.1572) (-2.6859,0.2365)
 (-2.6477,0.3052) (-2.6147,0.3648) (-2.6005,0.4281) (-2.5855,0.4849) (-2.5632,0.5308) (-2.5335,0.5674) (-2.4959,0.5971) (-2.4496,0.6228)
 (-2.3608,0.6632) (-2.2703,0.6981) (-2.1784,0.7282) (-2.0852,0.7536) (-1.9904,0.7743) (-1.8941,0.7895) (-1.7592,0.7998) (-1.6160,0.7970)
 (-1.4666,0.7948) (-1.3301,0.8426) (-1.1942,0.9007) (-1.0579,0.9643) (-0.9216,1.0328) (-0.7855,1.1050) (-0.6499,1.1797) (-0.5150,1.2555)
 (-0.3810,1.3307) (-0.2484,1.4037) (-0.1173,1.4727) (0.0117,1.5363) (0.1384,1.5930) (0.2624,1.6416) (0.3834,1.6813) (0.5017,1.7116)
 (0.6174,1.7323) (0.7314,1.7436) (0.8448,1.7458) (0.9595,1.7458) (1.0836,1.7824) (1.2141,1.8068) (1.3464,1.8108) (1.4438,1.8047)
 (1.5400,1.7931) (1.6353,1.7765) (1.7296,1.7551) (1.8229,1.7289) (1.9151,1.6971) (2.0184,1.6495) (2.1072,1.5841) (2.1773,1.5000)
 (2.2525,1.4344) (2.3378,1.3659) (2.4252,1.2749) (2.5131,1.1602) (2.5990,1.0229) (2.6805,0.8655) (2.7551,0.6919) (2.8209,0.5066)
 (2.8760,0.3146) (2.9189,0.1213) (2.9484,-0.0678) (2.9636,-0.2469) (2.9642,-0.4107) (2.9509,-0.5537) (2.9254,-0.6713) (2.8913,-0.7606)
 (2.8529,-0.8223) (2.8124,-0.8632) (2.7763,-0.9156) (2.7209,-0.9831) (2.6348,-1.0367) (2.5288,-1.0737) (2.4339,-1.0964) (2.3385,-1.1134)
 (2.2425,-1.1254) (2.1461,-1.1326) (2.0492,-1.1348) (1.9517,-1.1314) (1.8638,-1.1211) (1.7828,-1.0999) (1.7123,-1.0652) (1.6554,-1.0208)
 (1.6058,-0.9870) (1.5535,-0.9646) (1.4994,-0.9332) (1.4426,-0.8910) (1.3830,-0.8374) (1.3209,-0.7729) (1.2568,-0.6985) (1.1910,-0.6162)
 (1.1234,-0.5280) (1.0540,-0.4358) (0.9822,-0.3420) (0.9077,-0.2486) (0.8299,-0.1574) (0.7481,-0.0703) (0.6620,0.0111) (0.5714,0.0850)
 (0.4690,0.1377) (0.3553,0.1585) (0.2478,0.1712) (0.1485,0.1831) (0.0561,0.1967) (-0.0395,0.2162) (-0.1332,0.2411) (-0.2252,0.2708)
 (-0.3157,0.3050) (-0.4045,0.3439) (-0.4914,0.3883) (-0.5577,0.4260) (-0.6288,0.4658) (-0.7072,0.5047) (-0.7953,0.5384) (-0.8938,0.5603)
 (-0.9970,0.5566) (-1.0938,0.5216) (-1.1840,0.4736) (-1.2682,0.4169) (-1.3472,0.3539) (-1.4221,0.2863) (-1.4939,0.2161) (-1.5631,0.1449)
 (-1.6302,0.0744) (-1.6954,0.0064) (-1.7587,-0.0575) (-1.8200,-0.1158) (-1.8792,-0.1674) (-1.9360,-0.2116) (-1.9851,-0.2571) (-2.0243,-0.3156)
 (-2.0680,-0.3701) (-2.1193,-0.4157) (-2.1778,-0.4532)}
\begin{scope}[xshift=-4.9cm]
\path[fill=black!4,draw=black!20,line width=.3pt] (-3,-.85) -- (2.2,-.85) -- (3,.9) -- (-2.3,.9) -- cycle;
\draw[black!40,densely dashed,line width=.35pt] (1.65,.3) -- (1.65,1.85);
\draw[black!40,densely dashed,line width=.35pt] (2.25,-.3) -- (2.25,-1.2);
\fill (-2.5000,-0.6500) circle (.8pt);
\node at (-2.4800,-1.0500) {$p_1$};
\fill (-2.2000,0.8000) circle (.8pt);
\node at (-2.2400,1.1600) {$p_2$};
\fill (-0.2500,0.2000) circle (.8pt);
\node at (-0.2500,0.6000) {$p_3$};
\fill (1.6500,1.8500) circle (.8pt);
\node at (1.6700,2.2500) {$p_4$};
\fill (2.2500,-1.2000) circle (.8pt);
\node at (2.2700,-1.6000) {$p_5$};
\node at (0,-2.60) {(a) Prescribed points};
\end{scope}
\begin{scope}
\draw[black!75,line width=.55pt] plot[smooth cycle] coordinates {\constructioncurve};
\fill (-2.5000,-0.6500) circle (.8pt);
\node at (-2.4800,-1.0500) {$p_1$};
\fill (-2.2000,0.8000) circle (.8pt);
\node at (-2.2400,1.1600) {$p_2$};
\fill (-0.2500,0.2000) circle (.8pt);
\node at (-0.2500,0.6000) {$p_3$};
\fill (1.6500,1.8500) circle (.8pt);
\node at (1.6700,2.2500) {$p_4$};
\fill (2.2500,-1.2000) circle (.8pt);
\node at (2.2700,-1.6000) {$p_5$};
\node at (0,-2.60) {(b) A smooth closed curve};
\end{scope}
\begin{scope}[xshift=4.9cm]
\path[fill=black!9,even odd rule] plot[smooth cycle] coordinates {\constructionouter} plot[smooth cycle] coordinates {\constructioninner};
\draw[white,line width=1.6pt,opacity=.75] plot[smooth cycle] coordinates {\constructioncurve};
\draw[black!70,line width=.45pt] plot[smooth cycle] coordinates {\constructionouter};
\draw[black!70,line width=.45pt] plot[smooth cycle] coordinates {\constructioninner};
\draw[black!40,line width=.35pt] plot[smooth] coordinates {
 (-2.4818,-0.7228) (-2.4857,-0.7228) (-2.4899,-0.7210) (-2.4944,-0.7175) (-2.4990,-0.7123) (-2.5036,-0.7056) (-2.5081,-0.6975) (-2.5125,-0.6882)
 (-2.5165,-0.6780) (-2.5202,-0.6671) (-2.5233,-0.6558) (-2.5258,-0.6444) (-2.5278,-0.6331) (-2.5290,-0.6222) (-2.5295,-0.6119) (-2.5293,-0.6027)
 (-2.5284,-0.5946) (-2.5268,-0.5878) (-2.5245,-0.5826) (-2.5216,-0.5790) (-2.5182,-0.5772)};
\draw[black!30,densely dashed,line width=.3pt] plot[smooth] coordinates {
 (-2.5182,-0.5772) (-2.5143,-0.5772) (-2.5101,-0.5790) (-2.5056,-0.5825) (-2.5010,-0.5877) (-2.4964,-0.5944) (-2.4919,-0.6025) (-2.4875,-0.6118)
 (-2.4835,-0.6220) (-2.4798,-0.6329) (-2.4767,-0.6442) (-2.4742,-0.6556) (-2.4722,-0.6669) (-2.4710,-0.6778) (-2.4705,-0.6881) (-2.4707,-0.6973)
 (-2.4716,-0.7054) (-2.4732,-0.7122) (-2.4755,-0.7174) (-2.4784,-0.7210) (-2.4818,-0.7228)};
\draw[black!40,line width=.35pt] plot[smooth] coordinates {
 (-2.2216,0.8718) (-2.2177,0.8720) (-2.2134,0.8705) (-2.2088,0.8671) (-2.2039,0.8622) (-2.1990,0.8557) (-2.1941,0.8478) (-2.1893,0.8388)
 (-2.1848,0.8288) (-2.1807,0.8180) (-2.1770,0.8069) (-2.1739,0.7956) (-2.1715,0.7844) (-2.1697,0.7735) (-2.1687,0.7634) (-2.1685,0.7541)
 (-2.1691,0.7459) (-2.1704,0.7391) (-2.1724,0.7338) (-2.1751,0.7301) (-2.1784,0.7282)};
\draw[black!30,densely dashed,line width=.3pt] plot[smooth] coordinates {
 (-2.1784,0.7282) (-2.1823,0.7280) (-2.1866,0.7295) (-2.1912,0.7329) (-2.1961,0.7378) (-2.2010,0.7443) (-2.2059,0.7522) (-2.2107,0.7612)
 (-2.2152,0.7712) (-2.2193,0.7820) (-2.2230,0.7931) (-2.2261,0.8044) (-2.2285,0.8156) (-2.2303,0.8265) (-2.2313,0.8366) (-2.2315,0.8459)
 (-2.2309,0.8541) (-2.2296,0.8609) (-2.2276,0.8662) (-2.2249,0.8699) (-2.2216,0.8718)};
\draw[black!40,line width=.35pt] plot[smooth] coordinates {
 (1.6647,1.9235) (1.6682,1.9219) (1.6713,1.9185) (1.6738,1.9134) (1.6757,1.9067) (1.6770,1.8987) (1.6777,1.8894) (1.6776,1.8792)
 (1.6769,1.8682) (1.6755,1.8569) (1.6735,1.8453) (1.6709,1.8338) (1.6678,1.8228) (1.6643,1.8124) (1.6604,1.8030) (1.6562,1.7947)
 (1.6519,1.7877) (1.6476,1.7823) (1.6433,1.7786) (1.6392,1.7766) (1.6353,1.7765)};
\draw[black!30,densely dashed,line width=.3pt] plot[smooth] coordinates {
 (1.6353,1.7765) (1.6318,1.7781) (1.6287,1.7815) (1.6262,1.7866) (1.6243,1.7933) (1.6230,1.8013) (1.6223,1.8106) (1.6224,1.8208)
 (1.6231,1.8318) (1.6245,1.8431) (1.6265,1.8547) (1.6291,1.8662) (1.6322,1.8772) (1.6357,1.8876) (1.6396,1.8970) (1.6438,1.9053)
 (1.6481,1.9123) (1.6524,1.9177) (1.6567,1.9214) (1.6608,1.9234) (1.6647,1.9235)};
\draw[black!40,line width=.35pt] plot[smooth] coordinates {
 (2.2575,-1.2746) (2.2536,-1.2741) (2.2497,-1.2717) (2.2458,-1.2676) (2.2420,-1.2618) (2.2384,-1.2545) (2.2351,-1.2458) (2.2321,-1.2360)
 (2.2296,-1.2253) (2.2276,-1.2140) (2.2261,-1.2024) (2.2252,-1.1907) (2.2250,-1.1792) (2.2253,-1.1682) (2.2263,-1.1581) (2.2278,-1.1489)
 (2.2299,-1.1410) (2.2325,-1.1346) (2.2355,-1.1298) (2.2389,-1.1267) (2.2425,-1.1254)};
\draw[black!30,densely dashed,line width=.3pt] plot[smooth] coordinates {
 (2.2425,-1.1254) (2.2464,-1.1259) (2.2503,-1.1283) (2.2542,-1.1324) (2.2580,-1.1382) (2.2616,-1.1455) (2.2649,-1.1542) (2.2679,-1.1640)
 (2.2704,-1.1747) (2.2724,-1.1860) (2.2739,-1.1976) (2.2748,-1.2093) (2.2750,-1.2208) (2.2747,-1.2318) (2.2737,-1.2419) (2.2722,-1.2511)
 (2.2701,-1.2590) (2.2675,-1.2654) (2.2645,-1.2702) (2.2611,-1.2733) (2.2575,-1.2746)};
\draw[black!40,line width=.35pt] plot[smooth] coordinates {
 (-0.2748,0.1292) (-0.2780,0.1313) (-0.2806,0.1351) (-0.2824,0.1405) (-0.2834,0.1474) (-0.2835,0.1556) (-0.2829,0.1648) (-0.2814,0.1749)
 (-0.2792,0.1857) (-0.2762,0.1968) (-0.2727,0.2079) (-0.2685,0.2189) (-0.2639,0.2294) (-0.2589,0.2392) (-0.2538,0.2480) (-0.2485,0.2557)
 (-0.2433,0.2619) (-0.2382,0.2667) (-0.2334,0.2698) (-0.2291,0.2712) (-0.2252,0.2708)};
\draw[black!30,densely dashed,line width=.3pt] plot[smooth] coordinates {
 (-0.2252,0.2708) (-0.2220,0.2687) (-0.2194,0.2649) (-0.2176,0.2595) (-0.2166,0.2526) (-0.2165,0.2444) (-0.2171,0.2352) (-0.2186,0.2251)
 (-0.2208,0.2143) (-0.2238,0.2032) (-0.2273,0.1921) (-0.2315,0.1811) (-0.2361,0.1706) (-0.2411,0.1608) (-0.2462,0.1520) (-0.2515,0.1443)
 (-0.2567,0.1381) (-0.2618,0.1333) (-0.2666,0.1302) (-0.2709,0.1288) (-0.2748,0.1292)};
\draw[black!40,line width=.35pt] plot[smooth] coordinates {
 (-3.3138,-0.3306) (-3.3101,-0.3194) (-3.3008,-0.3088) (-3.2863,-0.2992) (-3.2667,-0.2908) (-3.2427,-0.2838) (-3.2148,-0.2784) (-3.1836,-0.2746)
 (-3.1500,-0.2727) (-3.1148,-0.2727) (-3.0789,-0.2745) (-3.0431,-0.2780) (-3.0083,-0.2834) (-2.9754,-0.2903) (-2.9451,-0.2986) (-2.9183,-0.3081)
 (-2.8956,-0.3186) (-2.8775,-0.3299) (-2.8646,-0.3415) (-2.8570,-0.3534) (-2.8550,-0.3651)};
\draw[black!30,densely dashed,line width=.3pt] plot[smooth] coordinates {
 (-2.8550,-0.3651) (-2.8587,-0.3763) (-2.8680,-0.3869) (-2.8825,-0.3965) (-2.9021,-0.4049) (-2.9261,-0.4119) (-2.9540,-0.4173) (-2.9852,-0.4211)
 (-3.0188,-0.4230) (-3.0540,-0.4230) (-3.0899,-0.4212) (-3.1257,-0.4176) (-3.1605,-0.4123) (-3.1934,-0.4054) (-3.2237,-0.3971) (-3.2505,-0.3876)
 (-3.2732,-0.3771) (-3.2913,-0.3658) (-3.3042,-0.3542) (-3.3118,-0.3423) (-3.3138,-0.3306)};
\draw[black!40,line width=.35pt] plot[smooth] coordinates {
 (-0.7070,1.6753) (-0.6956,1.6785) (-0.6817,1.6766) (-0.6656,1.6698) (-0.6478,1.6582) (-0.6286,1.6421) (-0.6086,1.6219) (-0.5883,1.5980)
 (-0.5681,1.5711) (-0.5485,1.5419) (-0.5301,1.5110) (-0.5132,1.4792) (-0.4983,1.4472) (-0.4858,1.4160) (-0.4760,1.3862) (-0.4691,1.3587)
 (-0.4652,1.3339) (-0.4646,1.3127) (-0.4671,1.2954) (-0.4727,1.2825) (-0.4814,1.2744)};
\draw[black!30,densely dashed,line width=.3pt] plot[smooth] coordinates {
 (-0.4814,1.2744) (-0.4928,1.2712) (-0.5067,1.2731) (-0.5228,1.2799) (-0.5406,1.2915) (-0.5598,1.3076) (-0.5798,1.3278) (-0.6001,1.3517)
 (-0.6203,1.3786) (-0.6399,1.4078) (-0.6583,1.4387) (-0.6752,1.4705) (-0.6901,1.5025) (-0.7026,1.5337) (-0.7124,1.5635) (-0.7193,1.5911)
 (-0.7232,1.6158) (-0.7238,1.6370) (-0.7213,1.6543) (-0.7157,1.6672) (-0.7070,1.6753)};
\draw[black!40,line width=.35pt] plot[smooth] coordinates {
 (3.2035,0.8072) (3.2049,0.7955) (3.2009,0.7820) (3.1917,0.7672) (3.1774,0.7514) (3.1586,0.7349) (3.1355,0.7183) (3.1087,0.7019)
 (3.0790,0.6861) (3.0471,0.6713) (3.0137,0.6579) (2.9797,0.6462) (2.9458,0.6365) (2.9131,0.6290) (2.8821,0.6239) (2.8538,0.6213)
 (2.8288,0.6213) (2.8076,0.6240) (2.7910,0.6291) (2.7792,0.6367) (2.7725,0.6465)};
\draw[black!30,densely dashed,line width=.3pt] plot[smooth] coordinates {
 (2.7725,0.6465) (2.7711,0.6583) (2.7751,0.6718) (2.7843,0.6866) (2.7985,0.7024) (2.8174,0.7188) (2.8405,0.7354) (2.8672,0.7518)
 (2.8969,0.7676) (2.9289,0.7824) (2.9623,0.7958) (2.9963,0.8076) (3.0301,0.8173) (3.0629,0.8248) (3.0938,0.8299) (3.1222,0.8325)
 (3.1472,0.8324) (3.1683,0.8298) (3.1850,0.8246) (3.1968,0.8170) (3.2035,0.8072)};
\draw[black!40,line width=.35pt] plot[smooth] coordinates {
 (0.9602,-1.0596) (0.9547,-1.0491) (0.9536,-1.0351) (0.9569,-1.0180) (0.9645,-0.9981) (0.9763,-0.9760) (0.9919,-0.9522) (1.0110,-0.9273)
 (1.0330,-0.9019) (1.0576,-0.8767) (1.0839,-0.8522) (1.1115,-0.8291) (1.1396,-0.8078) (1.1675,-0.7891) (1.1946,-0.7732) (1.2201,-0.7607)
 (1.2435,-0.7518) (1.2641,-0.7467) (1.2816,-0.7455) (1.2953,-0.7484) (1.3050,-0.7551)};
\draw[black!30,densely dashed,line width=.3pt] plot[smooth] coordinates {
 (1.3050,-0.7551) (1.3105,-0.7656) (1.3117,-0.7796) (1.3084,-0.7968) (1.3008,-0.8166) (1.2890,-0.8387) (1.2734,-0.8625) (1.2543,-0.8874)
 (1.2322,-0.9128) (1.2077,-0.9380) (1.1813,-0.9625) (1.1538,-0.9857) (1.1257,-1.0069) (1.0978,-1.0256) (1.0707,-1.0415) (1.0452,-1.0540)
 (1.0218,-1.0629) (1.0011,-1.0680) (0.9837,-1.0692) (0.9700,-1.0663) (0.9602,-1.0596)};
\draw[black!40,line width=.35pt] plot[smooth] coordinates {
 (-1.0626,-0.0082) (-1.0731,-0.0137) (-1.0871,-0.0147) (-1.1043,-0.0114) (-1.1241,-0.0037) (-1.1462,0.0081) (-1.1699,0.0237) (-1.1948,0.0429)
 (-1.2201,0.0650) (-1.2453,0.0896) (-1.2697,0.1160) (-1.2928,0.1436) (-1.3139,0.1718) (-1.3326,0.1997) (-1.3484,0.2268) (-1.3609,0.2524)
 (-1.3698,0.2758) (-1.3748,0.2964) (-1.3759,0.3139) (-1.3731,0.3276) (-1.3663,0.3373)};
\draw[black!30,densely dashed,line width=.3pt] plot[smooth] coordinates {
 (-1.3663,0.3373) (-1.3558,0.3428) (-1.3418,0.3439) (-1.3246,0.3406) (-1.3048,0.3329) (-1.2827,0.3211) (-1.2590,0.3054) (-1.2341,0.2863)
 (-1.2088,0.2642) (-1.1836,0.2396) (-1.1592,0.2132) (-1.1361,0.1855) (-1.1150,0.1574) (-1.0963,0.1294) (-1.0805,0.1023) (-1.0680,0.0768)
 (-1.0591,0.0534) (-1.0541,0.0327) (-1.0530,0.0153) (-1.0558,0.0015) (-1.0626,-0.0082)};
\draw[black!45,densely dashed,line width=.35pt] plot[smooth cycle] coordinates {\constructioncurve};
\fill (-2.5000,-0.6500) circle (.8pt);
\node at (-2.4800,-1.0500) {$p_1$};
\fill (-2.2000,0.8000) circle (.8pt);
\node at (-2.2400,1.1600) {$p_2$};
\fill (-0.2500,0.2000) circle (.8pt);
\node at (-0.2500,0.6000) {$p_3$};
\fill (1.6500,1.8500) circle (.8pt);
\node at (1.6700,2.2500) {$p_4$};
\fill (2.2500,-1.2000) circle (.8pt);
\node at (2.2700,-1.6000) {$p_5$};
\node at (0,-2.60) {(c) A tube with five necks};
\end{scope}
\end{tikzpicture}
\caption{A schematic view of the construction in $\R^3$.}
\label{fig:construction}
\end{figure}
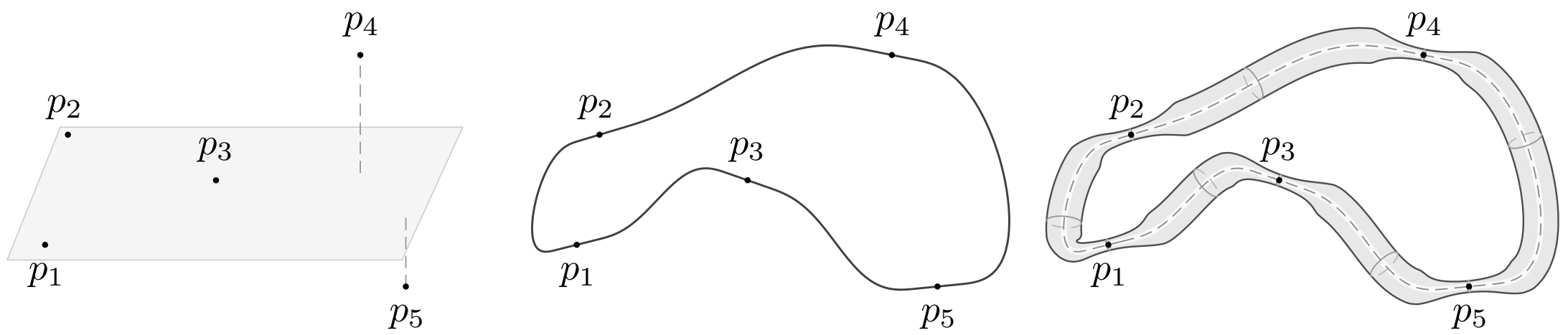

Choose a smooth closed embedded curve through the prescribed points, straight near each $p_i$. In a thin tube around this curve, insert a copy of the rotational neckpinch model of Proposition~\ref{prop:one-neck-reference} at scale $\sqrt T$ near each $p_i$. Its cylindrical ends join smoothly to the tube. Since joining the necks can shift where and when they pinch, we allow each neck to move along its straight segment, the segment to move transversely, and the neck scale to vary. This gives four parameters per neck (Figure~\ref{fig:construction}).

For a fixed flow, we recenter at each neck at successive steps. In the original coordinates these centers converge to the point where that neck pinches. The mean-convex level-set flow allows us to follow each neck even if another has already pinched. The centering estimates show that the differences between the resulting and prescribed singular points and times depend continuously on the initial parameters. In rescaled coordinates, convergence to the model flows and spacetime isolation of their neckpinches imply that this map converges uniformly to a positive multiple of the identity as $T\to0$. Degree theory then selects the initial parameters for which the necks pinch at the prescribed points and time, and uniform convergence forces these parameters to tend to zero.

For these parameters, compactness after rescaling by $T^{-1/2}$ forces any additional singularities at or before $T$ to approach one of the model neckpinches. We exclude them by proving that each prescribed neckpinch is the only singularity in a backward spacetime neighborhood whose size is independent of $T$ in these coordinates.

\subsection*{Acknowledgment}
The author is supported by the University of Macau Research Grant (MYRG-\allowbreak SRG2026-\allowbreak 00039-\allowbreak FSC). The author thanks Ao Sun for sharing this problem.

\subsection*{AI declaration}
Before using AI, the author had developed the main ideas and proof strategy and completed a draft. ChatGPT 6 Pro was used to assist with literature searches, streamlining the proof, and checking the arguments for possible gaps. The original proof used the Wa\.zewski box argument, and ChatGPT suggested replacing it with a direct degree argument. The author checked and revised the arguments and takes full responsibility for the mathematical content and final exposition.

\section{Preliminaries}\label{sec:preliminaries}
To continue mean curvature flow through singularities, Chen--Giga--Goto \cite{ChenGigaGoto1991} and Evans--Spruck \cite{EvansSpruck1991} independently developed the level-set formulation. Let $K_t$ be the level-set flow of a smooth strictly mean-convex domain $K_0$ and write $M_t=\partial K_t$. Evans--Spruck proved nonfattening for $\partial K_0$ when $K_0$ is compact. For compact mean-convex $K_0$, White \cite{White2000} proved that $\partial K_t$ is the level-set flow of $\partial K_0$ and $K_{t'}\subset\operatorname{Int}K_t$ for $t'>t$.

Colding--Minicozzi \cite{ColdingMinicozzi2012} define the entropy of a surface $M\subset\R^3$ by
\[
 \lambda(M)=\sup_{x_0\in\R^3,\,r>0}
 \frac{1}{4\pi r^2}\int_M e^{-|x-x_0|^2/(4r^2)}\dd\mathcal H^2(x).
\]
By Huisken \cite{Huisken1990}, entropy is nonincreasing, so $\lambda(M_0)$ gives uniform area-ratio bounds.

For $\alpha>0$, we call a smooth strictly mean-convex domain $K$ and its boundary $\alpha$-noncollapsed if at each $p\in\partial K$ there are interior and exterior tangent balls of radius at least $\alpha/H(p)$. Andrews \cite{Andrews2012} proved that smooth mean curvature flow of a closed mean-convex hypersurface preserves this condition with the same $\alpha$. Haslhofer--Kleiner \cite{HK} proved that the level-set flow from a smooth compact $\alpha$-noncollapsed domain retains the same bound through singularities. They also show that every sequence of such flows with fixed $\alpha$ has a subsequence whose spacetime domains and boundaries converge locally in the Hausdorff sense. The convergence is smooth near every positive-time regular boundary point of the limit. Thus any positive-time limit of singular points is singular.

Let $(p,T)$ be a spacetime point. For $t<T$, set $\tau=-\log(T-t)$ and consider the rescaled mean curvature flow $(T-t)^{-1/2}(M_t-p)$. When this flow is locally a normal graph over $\R\times\mathbb S^1(\sqrt2)$, we write its parametrization as
\[
\bigl(z,(\sqrt2+u(z,\theta,\tau))\cos\theta,(\sqrt2+u(z,\theta,\tau))\sin\theta\bigr).
\]

The graph equation of Zhou Gang, Knopf, and Sigal \cite{GangKnopfSigal2018} gives
\begin{equation}\label{eq:exact}
 u_\tau=Lu+\mathcal Q(u),
\end{equation}
where $L=\partial_z^2-\frac z2\partial_z+\frac12\partial_\theta^2+1$. The nonlinearity $\mathcal Q$ has no constant or linear terms and is affine in the second derivatives of $u$.

We use the Gaussian inner product and Sobolev norms
\[
 \left\langle f,g\right\rangle_G=\int_{\R\times\mathbb S^1}fg\,e^{-z^2/4}\dd z\dd\theta,
 \qquad
 \|f\|_{H_G^k}^2=\sum_{a+b\leq k}
       \|\partial_z^a\partial_\theta^b f\|_{L_G^2}^2.
\]

Let $E_-$ be the span of the eigenfunctions of $-L$ with negative eigenvalues. The spectral decomposition in \cite{SXI} gives
\[
 E_-=\operatorname{span}\{1,z,\cos\theta,\sin\theta\},\qquad
 \ker L=\operatorname{span}\{H_2,z\cos\theta,z\sin\theta\}.
\]
Here $H_2(z)=z^2-2$ is the quadratic Hermite polynomial. Set $E_+=(E_-\oplus\ker L)^\perp$; the spectrum of $-L$ on $E_+$ is contained in $[1/2,\infty)$. Let $P_-,P_+$ and $P_{\rm rot}$ be the Gaussian orthogonal projections onto $E_-,E_+$ and $\operatorname{span}\{z\cos\theta,z\sin\theta\}$, respectively. For $a\in E_-$, write $a=a_0+a_z z+a_\perp\cdot(\cos\theta,\sin\theta)$.

At a nondegenerate neckpinch in $\R^3$, Sun--Xue \cite[Thm.~1.2]{SXI} give the leading term $H_2/(2\sqrt2\,\tau)$ of $u$ in Gaussian $H^1$ on balls of radius proportional to $\sqrt\tau$. The quadratic part of $\mathcal Q$ evaluated at $u=bH_2$ is $-b^2H_2^2/(2\sqrt2)$, whose $H_2$ coefficient is $-\kappa b^2$ with $\kappa=2\sqrt2$. In their stability proof, Sun--Xue \cite{SXII} remove the $E_-$ component by translations and dilations at successive rescaled times. They denote this space by $E^+$, using the spectrum of $L$.

The projections are bounded in $H_G^k$. For every integer $k\geq0$ and $t>0$, expansion in eigenfunctions of $-L$ gives
\begin{equation}\label{eq:semigroup}
 \|e^{tL}P_+f\|_{H_G^{k+1}}
 \leq C(1+t^{-1/2})e^{-t/4}\|P_+f\|_{H_G^k},
\end{equation}
where $C$ depends only on $k$.

Ilmanen \cite{Ilmanen1994} showed that in the nonfattening case the evolving boundary is the spacetime support of the flow obtained by elliptic regularization. The estimate of Sun--Xue \cite[Thm.~2.4]{SXI} thus applies to our level-set flows in the following form.
\begin{lemma}\label{lem:cylindrical-pseudolocality}
Let $\{M_h\}_{0\leq h\leq h_0}$ be the rescaled boundary of an $\alpha$-noncollapsed mean-convex level-set flow with entropy at most $\Lambda$. Fix $0<h_1<h_0$, $0<\rho<1$, and a sufficiently small $\delta>0$. There are $r_0,\eta_0>0$, depending only on these parameters and $\alpha,\Lambda$, with the following property. Suppose that for some $r>r_0$ the surface $M_0$ agrees in $B_r(0)$ with the normal graph of $u(0)$ on $|z|<r$, and that $\|u(0)\|_{C^1(|z|<r)}\leq\eta_0$. For $0\leq h\leq h_0$, the surface $M_h$ agrees in $B_{\rho e^{h/2}r}(0)$ with a normal graph $u(h)$ on $|z|<\rho e^{h/2}r$, with
\[
 \|u(\cdot,h)\|_{C^3(|z|<\rho e^{h/2}r)}\leq\delta,
 \qquad h_1\leq h\leq h_0.
\]
\end{lemma}

We use the estimates of Sz{\'e}kelyhidi \cite[Props.~21--23]{SzekelyhidiSurgery2026} in the surface case.
\begin{lemma}\label{lem:neckpinch-comparison}
Fix a mean-convex flow $\Sigma_\tau$ in $\R^3$, rescaled about a nondegenerate neckpinch whose axis is the $z$-axis.
\begin{enumerate}
\item[(i)] For some $d>0$ and all sufficiently large $\tau$, the normal graph function of $\Sigma_\tau$ over the round cylinder on $B_{\tau^d}$ is $H_2/(2\sqrt2\,\tau)+O(\tau^{-19/10})$ in $C^2$. Translating $\Sigma_\tau$ by $-V\sqrt\tau$ along the axis gives, on each fixed ball, a graph over the round cylinder of radius $\sqrt{2+V^2}$, with graph function tending to zero in $C^1$ uniformly for $|V|\leq1$.
\item[(ii)] Suppose a rescaled mean-convex flow with entropy at most $\Lambda$ is a normal graph over $\Sigma_\tau$ on $B_{\tau^d}$ for some $d>0$ and all $\tau\geq T_0$, and its height has $L^2$ norm less than $\tau^{-9/8}$ with respect to $e^{-|X|^2/4}\dd\mathcal H^2(X)$ on $\Sigma_\tau\cap B_{\tau^d}$. Suppose also that on $B_{K\sqrt\tau}$ the flow is a normal graph over $\Sigma_\tau$ with $C^4$ norm at most $\delta$ for $T_0\leq\tau\leq3T_0$, where $0<K<2$. If $\delta$ is sufficiently small and $K\sqrt{T_0}$ sufficiently large, the graph extends to $B_{2\sqrt\tau}$ with the same $C^4$ bound for all $\tau\geq T_1$. Here $T_1$ depends only on $\Sigma_\tau$, $\Lambda$ and the parameters $d,\delta,T_0,K$.
\end{enumerate}
\end{lemma}
For each fixed $K_0>0$, Proposition~21(1) gives local cylindrical approximations to $\Sigma_\tau$ and the bound $|A_{\Sigma_\tau}|^2\leq 1/2+o(1)$ on $B_{K_0\sqrt\tau}$ as $\tau\to\infty$. Applying cylindrical pseudolocality on these patches and covering the ball as in the proof of Proposition~22 gives the extension estimate used in the proof of Proposition~23, which yields (ii) since $1/2<(9/8)/2$. When the graph is cut off on a larger ball in that proof, its bounded height contributes only an exponentially small Gaussian tail outside $B_{\tau^d}$.

\section{The rotational model and initial surfaces}\label{sec:global-family}
For a positive function $R$ on $\R$, let $M^R$ be the surface of revolution parametrized by $(z,R(z)\cos\theta,R(z)\sin\theta)$. We need a profile with cylindrical ends so that the necks join smoothly to the tube. Section~\ref{sec:local-estimates} uses the quadratic term in the model profile to determine the axial center of each perturbed neck.
\begin{proposition}
\label{prop:one-neck-reference}
There is a smooth positive even function $R$ on $\R$, equal to a constant $a_\infty>\sqrt2$ outside a bounded interval, such that $M^R$ is strictly mean-convex, has finite entropy, and is $\alpha$-noncollapsed for some $\alpha>0$. Its level-set flow $M_t^R$ has its unique first singularity at $(0,1)$. The radius of $(1-t)^{-1/2}M_t^R$ is $\sqrt2+u(z,\tau)$, where $\tau=-\log(1-t)$ and
\[
 u(z,\tau)=\frac{z^2-2}{2\sqrt2\,\tau}+o(\tau^{-1})
\]
uniformly on each compact axial interval as $\tau\to\infty$.
\end{proposition}

\begin{proof}
Sun--Xue \cite[(3.12)--(3.13)]{SXII} construct a nondegenerate neckpinch at $(0,S)$ from a smooth even profile $R_0$, nondecreasing for $z\geq0$ and constant for large $|z|$. Join the neck profile to the constant ends by an interpolation with uniformly bounded first and second derivatives in $z/\sqrt L$, where $L$ is the large parameter in their construction. Since $R_0=O(1)$ and $R_0''=O(L^{-1})$, we have $R_0R_0''<1+(R_0')^2$ for large $L$. The surface is therefore strictly mean-convex. Its cylindrical ends give finite entropy, $\alpha$-noncollapsing for some $\alpha>0$, and $\inf H>0$. Let $K_t$ be the level-set flow of the solid tube $K_0$ bounded by $M^{R_0}$. For small $h>0$, $K_h$ lies a positive distance inside $K_0$. By comparison as in White \cite[Cor.~3.3]{White2000}, the level-set flow from $\partial K_0$ at time $t$ is disjoint from $K_{t+h}$. Comparison with shrinking balls puts each interior point of $K_t$ in $K_{t+h}$ for some such $h$, proving nonfattening.

Let $R(z,t)$ be the evolving profile with initial value $R_0(z)$. At positive times before the first singularity, the strong maximum principle gives $R_z>0$ on $z>0$ and, by reflection, $R_z<0$ on $z<0$. The local arguments of Altschuler--Angenent--Giga \cite[Lemmas~5.2--5.3]{AltschulerAngenentGiga1995} give a positive lower bound for the radius on each compact subinterval of either half-line and smoothness wherever the limiting radius is positive. Thus the first singularity occurs when the central circle collapses to the origin. An earlier collapse would leave at most the origin in the plane $z=0$ thereafter, since the domains decrease with time. The local neck constructed by Sun--Xue has positive radius at $z=0$ for $t$ just below $S$, ruling this out.

Since the radius at $z\ne0$ stays positive at $S$, the enclosing shrinking cylinder survives past $S$ by comparison. Parabolic dilation by $S^{-1/2}$ puts the first singularity at $(0,1)$ and gives $R(z)=S^{-1/2}R_0(\sqrt S\,z)$ with $a_\infty>\sqrt2$. Sun--Xue \cite[Cor.~1.4]{SXI} give the quadratic term.
\end{proof}

Choose disjoint balls about $p_i$ and a diameter through $p_i$ in each. In $\R^3$ their endpoints can be joined cyclically by mutually disjoint smooth arcs outside the balls. Smoothing the joins outside the straight segments gives a smooth embedded closed curve $\Gamma$ through all $p_i$. Fix $\rho>0$ so that $\Gamma$ is straight for arclength $|z|<10\rho$ about each $p_i$ and the balls $B_{4\rho}(p_i)$ are disjoint and meet $\Gamma$ only along these straight segments. The normal injectivity radius of the compact curve $\Gamma$ is positive.

In the formulas for the $i$-th neck, we identify the model's $z$-axis with the corresponding straight segment of $\Gamma$. Fix $0<\eta<1/2$. Write the parameters as $\xi=(\zeta_i,\vartheta_i)_{i=1}^n\in Q=[-1,1]^{4n}$, with $\zeta_i\in[-1,1]^3$ and $\vartheta_i\in[-1,1]$, and set $\sigma_i=1+\eta\vartheta_i$. The corresponding translated and scaled model flows are
\[
 M_{i,t}^{R,\xi}
 =\eta\zeta_i+\sqrt{\sigma_i}\,
   M_{t/\sigma_i}^R.
\]
Each model flow has its unique first singularity at $(\eta\zeta_i,\sigma_i)$.

\begin{lemma}
\label{lem:global-tube-family}
For every sufficiently small $T>0$ there is a smooth family of embedded strictly mean-convex tori $M_0^{T,\xi}=\partial K_0^{T,\xi}$, $\xi\in[-1,1]^{4n}$, bounding tubes of radius at most $C\sqrt T$ about curves $\Gamma_{T,\xi}$. In $B_{\rho/2}(p_i)$, each torus agrees with $p_i+\sqrt T\,M_{i,0}^{R,\xi}$. The domains are $\alpha$-noncollapsed and their boundaries satisfy
\[
 |A|\leq CT^{-1/2},\qquad \lambda(M_0^{T,\xi})\leq\Lambda,
\]
with $\alpha,C,\Lambda>0$ independent of $T,\xi,\eta$.
\end{lemma}

\begin{proof}
Translate the straight piece near $p_i$ by $\sqrt T\,\eta\zeta_i$, using a smooth ambient vector field that is constant on $B_{3\rho}(p_i)$ and vanishes outside $B_{4\rho}(p_i)$. This produces a closed curve $\Gamma_{T,\xi}$ which converges to $\Gamma$ smoothly, uniformly in $\xi$, and whose normal injectivity radius is bounded below independently of $T$ and $\xi$.  On each translated straight piece, use its arclength coordinate $z$, centered at $p_i+\sqrt T\,\eta\zeta_i$.

Choose disjoint smooth arclength cutoffs $\chi_i$, equal to one on $|z|\leq\rho$ and zero on $|z|\geq2\rho$.  On the center curve prescribe the positive tube radius
\begin{equation}
\label{eq:global-radius-definition}
 \sqrt T\left[a_\infty+\sum_i\chi_i
 \left(\sqrt{\sigma_i}\,R\left(\frac{z}{\sqrt{T\sigma_i}}\right)-a_\infty\right)\right].
\end{equation}
The supports are disjoint and $R$ is bounded, so the tube radius is at most $C\sqrt T$. For small $T$, the radius on $\rho\leq|z|\leq2\rho$ is $\sqrt T\,a_\infty[1+\chi_i(\sqrt{\sigma_i}-1)]$, with derivatives $O(\sqrt T)$. The boundary of the normal disk bundle with radius \eqref{eq:global-radius-definition} is therefore a smooth embedded torus.

After rescaling by $T^{-1/2}$ about points of $\Gamma_{T,\xi}$, the center curves converge locally to straight lines, and the derivatives of the rescaled tube radius tend to zero where the cutoffs vary. After passing to a subsequence, the pointed initial domains converge smoothly either to a translate of the region bounded by $\sqrt\sigma M^R$, $\sigma\in[1-\eta,1+\eta]$, or, when the normalized distance to every neck tends to infinity, to a solid cylinder whose radius is bounded above and away from zero uniformly in $\eta,\xi$ by \eqref{eq:global-radius-definition}. The normal injectivity radius excludes additional sheets. These limits have uniformly positive mean curvature, bounded curvature derivatives and interior and exterior tangent balls of a fixed radius. Taking balls of half this radius gives strict separation from the surface to second order at the point of tangency and positive separation away from it, so the balls persist under smooth local convergence. Smooth convergence bounds $\sqrt T\,|A|$ above and $\sqrt T\,H$ away from zero. The tangent balls of radius $c\sqrt T$ thus give $\alpha$-noncollapsing uniformly in $T,\xi,\eta$.

The uniform normal injectivity radius and bounded total length give $\mathcal H^1(\Gamma_{T,\xi}\cap B_s(x))\leq Cs$ for $s\geq c\sqrt T$. Consequently,
\[
 \mathcal H^2(M_0^{T,\xi}\cap B_s(x))
 \leq C\sqrt T(s+\sqrt T)\leq Cs^2.
\]
For $s\leq c\sqrt T$ the same bound follows from the uniform normalized local graphs.  Summing the Gaussian integral over dyadic annuli proves $\lambda(M_0^{T,\xi})\leq\Lambda$.
\end{proof}

\begin{lemma}
\label{lem:canonical-pointed-convergence}
Write $M_t^{T,\xi}=\partial K_t^{T,\xi}$ for the boundary of the level-set flow from $K_0^{T,\xi}$. If $T_j\downarrow0$ and $\xi_j\to\xi$, then for every fixed $i$ the flows
\[
 T_j^{-1/2}
 \bigl(M_{T_j t}^{T_j,\xi_j}-p_i\bigr)
\]
converge to $M_{i,t}^{R,\xi}$. If $\xi_j\to0$ and $q_j\in\Gamma_{T_j,\xi_j}$ satisfy $\min_i|q_j-p_i|/\sqrt{T_j}\to\infty$, rescaling about $q_j$ instead gives, along a subsequence, a shrinking round cylinder of initial radius $a_\infty$. In both cases the spacetime domains and their boundaries converge locally in the Hausdorff sense, and the convergence is smooth near every positive-time regular point of the limiting flow.
\end{lemma}

\begin{proof}
At $p_i$, Lemma~\ref{lem:global-tube-family} gives smooth local convergence of the initial domains to the translated and scaled model. At $q_j$, every fixed rescaled ball eventually misses the nonconstant part of each neck profile. Since $\xi_j\to0$, all $\sigma_i\to1$, so \eqref{eq:global-radius-definition} gives limiting radius $a_\infty$ even where the cutoffs vary. The rescaled center curves converge subsequentially to a straight line. In both cases the truncated signed distance functions converge locally uniformly. Level-set stability and comparison identify the limits as the model flow, which is nonfattening by the proof of Proposition~\ref{prop:one-neck-reference}, and the shrinking cylinder, respectively. Section~\ref{sec:preliminaries} gives the stated convergence. The model inherits the noncollapsing and entropy bounds.
\end{proof}

\section{Local graphical estimates and centering}\label{sec:local-estimates}
\subsection{Weighted estimates and centering}
We choose the cylinder axis at each rescaled time so that the cutoff graph function has zero projection onto $\operatorname{span}\{z\cos\theta,z\sin\theta\}$. At each integer time we also translate and dilate to set its projection onto $E_-$ to zero. We call this choice of coordinates centering. We count one unit of rescaled time in the frame chosen at the start of each step. If step $j$ starts at physical time $t_j$ with spatial scale $\sqrt{s_j}$, iteration time $j+h$ corresponds to physical time $t_j+s_j(1-e^{-h})$ for $0\leq h\leq1$.

Fix $\gamma=\frac34$ and $\sigma=\frac14$, and choose $0<\beta<1/8$ sufficiently small.

\begin{lemma}\label{lem:smoothing}
For large $j$, let $u$ be the normal graph function of a rescaled mean curvature flow over the cylinder on $|z|<16j^\beta$ for an interval of length at most two, with $\|u\|_{C^3}\leq C_0j^{-\gamma}$. Let $\chi_s(z)=\chi(z/s^\beta)$, where $\chi$ is smooth, even, equal to one on $[-14,14]$ and supported in $(-15,15)$. Put $v=\chi_j u$ and choose the axis so that $P_{\rm rot}v=0$. We transport the frame defining $\theta$ so that its angular velocity has zero component along the moving axis. If $\|v(s_0)\|_{H_G^1}\leq C_1/j$ at the initial time $s_0$, then for $s\geq s_0+1/4$,
\[
 \|v(s)\|_{H_G^2}\leq Cj^{-1},
 \qquad
 \|\mathcal Q(v(s))\|_{L_G^2}
       \leq Cj^{-1-\gamma}.
\]
If $\|v(s_0)\|_{H_G^2}\leq C_2/j$, both estimates hold throughout the interval. The same conclusions hold for $v=\chi_s u$ when $|s-j|\leq2$. The constants depend only on the cutoff, $C_0$ and $C_1$ for $s\geq s_0+1/4$, and also on $C_2$ for earlier times.
\end{lemma}
\begin{proof}
For small $|u|+|Du|$, the nonlinearity has the form $\mathcal Q(u)=A^{ab}(u,Du)D_{ab}u+O(|u|^2+|Du|^2)$, where $A^{ab}=O(|u|+|Du|^2)$ and the lower-order term is smooth in $(u,Du)$. Cutting off and replacing $u$ by $v$ in the coefficients introduces an error $E$ on the annulus where the cutoff varies. For either cutoff this annulus lies at $|z|\asymp j^\beta$, and
\begin{equation}\label{eq:tails}
 \|E\|_{H_G^k}\leq C e^{-c j^{2\beta}},\qquad k=0,1.
\end{equation}
The cutoff and rotation terms are given in \cite[(2.2)--(2.4)]{SXII}. On the cutoff annulus, $E$ and $DE$ involve derivatives of $u$ of order at most three, including for the moving cutoff $\chi_s$. The Gaussian weight gives \eqref{eq:tails}. Projecting the equation onto the two rotation modes and differentiating $P_{\rm rot}v=0$ bounds the angular velocity $\Omega\in\R^2$ of the axis by $|\Omega|\leq C\|\mathcal Q(v)\|_{L_G^2}+Ce^{-c j^{2\beta}}$.

Since $D\mathcal Q(v)$ involves at most third spatial derivatives of $v$, the $C^3$ bound and the quadratic vanishing of $\mathcal Q$ give $\|\mathcal Q(v)\|_{H_G^1}+|\Omega|\leq Cj^{-2\gamma}$. The rotation generators and their first derivatives are bounded by $C(1+|z|)$ on the graph. Together with \eqref{eq:tails}, this gives $\|v_s-Lv\|_{H_G^1}\leq Cj^{-2\gamma}$.

On intervals of length $h\leq2$, Duhamel's formula and \eqref{eq:semigroup} bound the $H_G^2$ norm of the forcing contribution by $Cj^{-2\gamma}\int_0^h(1+t^{-1/2})\,dt$. The finitely many modes in $E_-\oplus\ker L$ obey the same bound. For $h\geq1/4$ the same estimate bounds the homogeneous term in $H_G^2$ by its initial $H_G^1$ norm times a fixed constant. For $h\leq1/4$ the Hermite--Fourier expansion bounds the $H_G^2$ operator norm by $e^h$. These estimates give the asserted smoothing and propagation since $2\gamma>1$. The quasilinear formula also gives $\|\mathcal Q(v)\|_{L_G^2}\leq C\|v\|_{C^1}\|v\|_{H_G^2}\leq Cj^{-1-\gamma}$.
\end{proof}

Fix $B>0$ for the remainder estimates below.
\begin{lemma}\label{lem:weighted-IFT}
For sufficiently large $j$, suppose $u$ defines a normal graph on $|z|<15j^\beta+1$, $\|u\|_{C^3}\le Cj^{-\gamma}$, and $v=\chi_j u$ for the cutoff in Lemma~\ref{lem:smoothing} satisfies
\[
 v=bH_2+a+w,\qquad c/j\le b\le C/j,\qquad
 \|w\|_{H_G^1}\le Bj^{-3/2},
\]
where $a\in E_-$, $w\in E_+$, and
\[
 |a_0|\le Cj^{-2},\qquad
 |a_z|+|a_{\perp}|\le Cj^{-2-\sigma}.
\]
A unique change of coordinates near the identity, given by a dilation $\lambda$, axial translation $y$, transverse translation $\mathbf{x}\in\R^2$ and axis rotation $\omega\in\R^2$, makes the new graph function $\tilde u$ satisfy $P_-(\chi_j\tilde u)=P_{\rm rot}(\chi_j\tilde u)=0$. Its parameters depend continuously on $u$ and satisfy
\begin{equation}\label{eq:parameters}
 |\lambda-1|\le Cj^{-2},\quad
 |y|\le Cj^{-1-\sigma},\quad
 |\mathbf{x}|+|\omega|\le Cj^{-2-\sigma}.
\end{equation}
\end{lemma}
\begin{proof}
Write $X=(z,Y)$ with $Y\in\R^2$ and use $X\mapsto\lambda Q_\omega(X-(y,\mathbf{x}))$, where $Q_\omega$ is generated by $(z,Y)\mapsto(-\omega\cdot Y,z\omega)$. Use the parameters $\mathbf p=(\lambda-1,\mathbf{x},by,\omega)$. Let $\Psi(u,\mathbf p)$ be the six coefficients to be removed after regraphing and applying $\chi_j$. The graph transformation formulas in \cite[Cor.~A.3]{SXII} give the leading variations
\[
 \sqrt2(\lambda-1),\quad -\mathbf{x}\cdot(\cos\theta,\sin\theta),
 \quad 2byz,\quad z\omega\cdot(\cos\theta,\sin\theta).
\]
The term $2byz$ shows that using $by$ as the axial parameter keeps the derivative of the $z$ coefficient bounded away from zero as $b\to0$. The dilation term $-b(z^2+2)$ is orthogonal to the five nonconstant modes appearing in $\Psi$, and the remaining rotation terms are $O(b)$ in Gaussian norms. Thus $D_{\mathbf p}\Psi$ is a fixed invertible matrix $D$ plus a small error.

At the identity, differentiating the graph function with respect to $by$ gives $u_z/b$. After cutoff, its difference from $2z$ is $(w_z+a_z)/b$ up to cutoff errors. The error in the $by$ column of $D_{\mathbf p}\Psi$ is at most
\[
 C\|w\|_{H_G^1}/b+C\|a\|_{L_G^2}/b+Ce^{-c j^{2\beta}}/b
 =O(j^{-1/2}).
\]
For $|\mathbf p|\leq C_1j^{-2}$, one has $|y|=O(j^{-1})$ and all other corrections are $O(j^{-2})$. The graph remains defined on $|z|<15j^\beta$, and the Gaussian measure changes by a factor whose value and reciprocal are at most
\[
 e^{C(|y|j^\beta+|\mathbf{x}|+|\lambda-1|j^{2\beta}+|\omega|j^{2\beta})}=1+o(1).
\]
On the cutoff support the inverse coordinate change displaces $(z,\theta)$ by at most $|y|+Cj^{-2}(1+|z|)=O(j^{-1})$. Thus the change in $u_z/b$ from evaluating it at the new point is bounded by
\[
 \frac{|y|+Cj^{-2}(1+|z|)}{b}\,\|u\|_{C^2}=O(j^{-\gamma}).
\]
Taking the normal component of the translation velocity $-\lambda Q_\omega e_z$, where $e_z$ is the axial unit vector, gives the $by$ derivative as $u_z/b$ at that point times $1+O(j^{-2}(1+|z|))$. Since $|u_z|/b=O(j^{1-\gamma})$, this factor changes the column by $o(j^{-\gamma})$. Cutoff errors remain exponentially small after division by $b$. In the other columns, the first variations involve only $u$ and $Du$. Their errors are $O(j^{-\gamma})$ after integration against the fixed eigenfunctions, since every polynomial in $z$ has finite Gaussian norm. Uniformly on this parameter ball,
\[
 \|D^{-1}(D_{\mathbf p}\Psi-D)\|
 \le C(j^{-1/2}+j^{-\gamma})+o(1)<1/2.
\]
Since $\Psi(u,0)=O(j^{-2})$, the contraction $\mathbf p\mapsto\mathbf p-D^{-1}\Psi(u,\mathbf p)$ preserves the ball for large $C_1$. Its unique fixed point is the desired zero and depends continuously on $u$.

In the five nonconstant rows, the dilation derivative is $O(j^{-\gamma})$, since its leading terms have zero projection. Inverting the remaining submatrix gives $|\mathbf{x}|+b|y|+|\omega|\leq Cj^{-2-\sigma}+Cj^{-\gamma}j^{-2}\leq Cj^{-2-\sigma}$, proving the sharper bounds.
\end{proof}

\begin{lemma}\label{lem:center}
At a sufficiently large integer time $j$, suppose $v=\chi_j u$ satisfies
\begin{equation}\label{eq:bootstrap}
 v=bH_2+w,\quad
 \frac{c_-}{j}\leq b\leq\frac{c_+}{j},\quad
 \|w\|_{H_G^1}\leq B j^{-3/2},\quad
 \|v\|_{H_G^2}\leq C_2/j,
\end{equation}
where $w\in E_+$, $0<c_-<\kappa^{-1}<c_+$ and $C_2$ is sufficiently large. On $[j,j+1]$, assume the $C^3$ bound of Lemma~\ref{lem:smoothing} and choose the axis so that $P_{\rm rot}v=0$. Its angular velocity satisfies
\begin{equation}\label{eq:omega}
 |\Omega(j+h)|\leq Cj^{-2-\sigma},\qquad 0\leq h\leq1.
\end{equation}
Before centering, $b(j+h)=b_j+O(j^{-2})$ throughout the step. At $j+1$ the centering satisfies \eqref{eq:parameters}, and the centered graph satisfies \eqref{eq:bootstrap} with $j$ replaced by $j+1$ and strict inequalities for $b$ and $\|w\|_{H_G^1}$. Its quadratic coefficient satisfies
\begin{equation}\label{eq:discrete-b}
 b_{j+1}=b_j-\kappa b_j^2+O(j^{-2-\sigma}).
\end{equation}
\end{lemma}
\begin{proof}
Write $v=b(s)H_2+a(s)+w(s)$, with $a(s)\in E_-$ and $w(s)\in E_+$. Work on the maximal subinterval of $[j,j+1]$ starting at $j$ on which $b\geq c_-/(2j)$, the upper bounds for $b$ and $\|w\|_{H_G^1}$ hold with twice their stated constants, and $\|a\|_{L_G^2}\leq2B j^{-3/2}$. Lemma~\ref{lem:smoothing} propagates the initial $H_G^2$ bound and gives $\|v\|_{H_G^2}\leq C/j$ and $\|\mathcal Q(v)\|_{L_G^2}\leq Cj^{-1-\gamma}$. On the round cylinder the rotational correction is a linear combination of $z\cos\theta$ and $z\sin\theta$, so its projection onto $E_+$ vanishes up to an exponentially small cutoff error. The remaining term has $L_G^2$ norm at most $C|\Omega|\|v\|_{H_G^2}$, using $\|zf\|_{L_G^2}\leq C(\|f_z\|_{L_G^2}+\|f\|_{L_G^2})$ from Gaussian integration by parts. This is $O(j^{-5/2})$ since $|\Omega|\leq Cj^{-2\gamma}$ by Lemma~\ref{lem:smoothing}. The Hermite--Fourier modes are orthogonal in $H_G^1$, and $-L$ has spectrum in $[1/2,\infty)$ on $E_+$. Hence the homogeneous evolution of $w(j)$ decays by a factor at most $e^{-h/2}$. Applying \eqref{eq:semigroup} to the projected forcing in Duhamel's formula gives, for $j+h$ in this interval,
\begin{equation}\label{eq:stable-improve}
 \|w(j+h)\|_{H_G^1}\leq e^{-h/2}\|w(j)\|_{H_G^1}
       +C\sqrt h\,j^{-1-\gamma}.
\end{equation}
Here $h=1$ denotes the value before centering. For large $j$ the right-hand side is less than $2Bj^{-3/2}$ throughout the interval.

Set $B_0=b\chi_jH_2$. Since $\mathcal Q(B_0)$ is even in $z$ and independent of $\theta$, only $\mathcal Q(v)-\mathcal Q(B_0)$ contributes to the projections of $\mathcal Q(v)$ onto $z,\cos\theta,\sin\theta$ and the two rotational modes. In $\mathcal Q(v)-\mathcal Q(B_0)$ the coefficients of the second derivatives of $v-B_0$ are $A^{ab}(v,Dv)$. These coefficients and their first derivatives are $O(j^{-\gamma})$, and the coefficients of the lower derivatives of $v-B_0$ are also $O(j^{-\gamma})$. Integrating the second derivatives by parts once against the Gaussian weight gives
\begin{equation}\label{eq:finite-lip}
 |\langle\mathcal Q(v)-\mathcal Q(B_0),\phi\rangle_G|
 \leq Cj^{-\gamma}\|v-B_0\|_{H_G^1}
       +Ce^{-c j^{2\beta}},
\end{equation}
where $\phi$ is one of $1,z,\cos\theta,\sin\theta,H_2,z\cos\theta,z\sin\theta$. Here $\phi$, $D\phi$ and $z\phi$ have finite Gaussian norms. Taylor expansion through second order and Gaussian integration, with cutoff errors bounded by \eqref{eq:tails}, give
\begin{equation}\label{eq:pure-moments}
\begin{aligned}
 \langle\mathcal Q(B_0),1\rangle_G&=O(j^{-2}),\\
 \langle\mathcal Q(B_0),H_2\rangle_G&=-\kappa b^2\|H_2\|_{L_G^2}^2+O(j^{-3}).
\end{aligned}
\end{equation}

Differentiating $P_{\rm rot}v=0$ gives a linear system for $\Omega$ whose matrix is the identity plus $O(j^{-1})$.  Its right-hand side is controlled by \eqref{eq:finite-lip}; hence
\[
 |\Omega|\leq Cj^{-\gamma}
          (\|a\|_{L_G^2}+\|w\|_{H_G^1})+Ce^{-c j^{2\beta}}.
\]
Since $a(j)=0$, variation of constants in the graph equation over one step gives
\[
 |a_0|\leq Cj^{-2},\qquad
 |a_z|+|a_\perp|+|\Omega|
       \leq Cj^{-3/2-\gamma}=Cj^{-2-\sigma}.
\]
This proves \eqref{eq:omega} and the endpoint bounds required by Lemma~\ref{lem:weighted-IFT}. Projection onto $H_2$, using \eqref{eq:finite-lip} and \eqref{eq:pure-moments}, gives $b'=-\kappa b^2+O(j^{-2-\sigma})$. Hence $b(j+h)=b(j)+O(j^{-2})$ before centering. These improved bounds for $a,b,w$ extend the interval to $j+1$.

At the endpoint, $15(j+1)^\beta+1<16j^\beta$ for large $j$. Replacing $\chi_j$ by $\chi_{j+1}$ changes the Gaussian projections by $O(e^{-c j^{2\beta}})$. The two rotation generators have invertible projection onto the rotational modes, so an axis rotation of the same order restores $P_{\rm rot}(\chi_{j+1}u)=0$. Applying Lemma~\ref{lem:weighted-IFT} with this cutoff gives \eqref{eq:parameters}.

The centering changes $v$ in $H_G^1$ by $O(j^{-2})$. Indeed, \eqref{eq:parameters} and the $H_G^2$ bound control translations; dilation contributes its constant radial term and $O(j^{-3+\beta})$ from $zv_z$; rotation contributes $O(j^{-2-\sigma})+O(j^{-3-\sigma+\beta})$. Here $\|zv_z\|_{H_G^1}\leq Cj^\beta\|v\|_{H_G^2}$ on the cutoff support, with exponentially small errors. Thus \eqref{eq:stable-improve} at $h=1$ strictly improves the bound on $\|w\|_{H_G^1}$. For the $H_G^2$ norm, second derivatives in the graph transformation appear linearly, lower derivatives are controlled by the small $C^3$ bound, and the Gaussian measure changes by a factor $1+o(1)$. The transformed norm is therefore at most $(1+o(1))\|v\|_{H_G^2}+O(j^{-2})$. The smoothing estimate from the initial $H_G^1$ bound is independent of $C_2$, so choosing $C_2$ larger than twice its constant gives $\|v_{j+1}\|_{H_G^2}\leq C_2/(j+1)$ for large $j$.

A dilation changes the quadratic coefficient by $O(b|\lambda-1|)=O(j^{-3})$; the constant radial term is orthogonal to $H_2$.  Translation of $bH_2$ has no first-order $H_2$ component.  On the remainder it contributes at most $C|y|\|w\|_{H_G^1}=O(j^{-5/2-\sigma})$.  The leading transverse and rotational generators are orthogonal to $H_2$ and their remaining terms are smaller still.  Hence the centering changes $b$ by $O(j^{-2-\sigma})$, giving \eqref{eq:discrete-b}. The map $b\mapsto b-\kappa b^2$ is increasing on $[c_-/j,c_+/j]$ for large $j$. At $b=c_\pm/j$ the right-hand side of \eqref{eq:discrete-b} differs from $c_\pm/(j+1)$ by $c_\pm(1-\kappa c_\pm)j^{-2}+O(j^{-2-\sigma})$, positive at the lower endpoint and negative at the upper. Thus $c_-/(j+1)<b_{j+1}<c_+/(j+1)$.

By \eqref{eq:parameters}, translations, dilation and rotation reduce the axial graphical radius by at most
\begin{equation}\label{eq:radius-loss}
 C\bigl(|y|+|\mathbf{x}|+j^\beta(|\lambda-1|+|\omega|)\bigr)\leq C'j^{-1-\sigma}.
\end{equation}
\end{proof}

\subsection{Graphical continuation and initialization}
To maintain the $C^3$ bound in Lemma~\ref{lem:center}, we first estimate the spatial derivatives in a fixed rescaled frame.
\begin{lemma}
\label{lem:propagate}
Let $u$ solve \eqref{eq:exact} in a fixed rescaled frame, and put $U(\xi,\theta,h)=u(e^{h/2}\xi,\theta,h)$. If the spatial norm $\|U(\cdot,h)\|_{C^3(|\xi|<R)}$ is small uniformly for $0\leq h\leq H$, where $R\geq1$ and $1\leq H\leq4$, then for $0<\rho<1$,
\[
 \|U(H)\|_{C^4(|\xi|\leq\rho R)}
       \leq C\|U(0)\|_{C^0(|\xi|<R)}+Ce^{-cR^2},
\]
with constants depending only on $\rho$ and the $C^3$ bound.
\end{lemma}
\begin{proof}
The change of variables cancels the axial drift. Since $H\leq4$, the equation for $U$ is uniformly parabolic. Absorb the gradient-square term into a bounded first-order coefficient and write the reaction term as a bounded coefficient times $U$. For small $\alpha>0$ and large $M$, the function $e^{Mh}[\|U(0)\|_\infty+\sum_{\pm}\exp\{-\alpha(R\pm\xi)^2/(h+1)\}]$ is a supersolution for $|U|$ dominating the initial and lateral data. On $|\xi|\leq(1+\rho)R/2$ comparison therefore bounds $|U|$ by $C\|U(0)\|_\infty+Ce^{-cR^2}$. Since $H\geq1$ and $(1-\rho)R/2\geq(1-\rho)/2$, quasilinear interior regularity bounds the coefficient derivatives on smaller positive-time cylinders. The interior estimates of Lieberman \cite{Lieberman1996} then give the stated $C^4$ bound.
\end{proof}

\begin{lemma}
\label{lem:extend}
Consider $\alpha$-noncollapsed level-set flows with entropy at most $\Lambda$. In the fixed frame before centering at $J$, assume that for $J-\log J\leq s\leq J$ the rescaled boundary agrees in $B_{32s^\beta}(0)$ with a normal graph on $|z|<32s^\beta$ satisfying $\|u(s)\|_{H_G^1}\leq C/s$ and $\|u(s)\|_{C^3}<C_0s^{-\gamma}$ there. Assume the centering at $J$ satisfies \eqref{eq:parameters}. Let $j\geq J$ be an integer and suppose the centered graph satisfies \eqref{eq:bootstrap} at each integer from $J$ through $j$. On each completed step $[k,k+1]$, $J\leq k<j$, assume the rescaled boundary agrees in $B_{16s^\beta}(0)$ with a normal graph on $|z|<16s^\beta$ satisfying $\|u(s)\|_{C^3}<C_0s^{-\gamma}$. For $J$ sufficiently large in terms of the common bounds, the same graphical description and $C^3$ bound hold on $[j,j+1]$.
\end{lemma}
\begin{proof}
As long as the graphical $C^3$ bound holds on $[j,s]\subset[j,j+1]$, Lemma~\ref{lem:smoothing} propagates the initial $H_G^2$ bound in \eqref{eq:bootstrap}, giving $\|\chi_j u(s)\|_{H_G^1}\leq C/s$. Replacing $\chi_j$ by $\chi_s$ costs $O(e^{-c j^{2\beta}})$ in $H_G^1$. Put $R_s=s^\beta$ and $Z=|\chi_s u|+|\nabla(\chi_s u)|$, so that $\|Z\|_{L_G^2}\leq C/s$. In a fixed frame, Sun--Xue \cite[Lem.~4.2]{SXI} give a differential inequality for $Z$ with cutoff error supported on $|z|\geq14R_s$. Fix the terminal time $s$ and write $r$ for time on $[s-h,s]$, where $h=2\log(7R_s)+2$. On this interval the moving axis adds a source bounded by $Cs^{-3/2}(1+z^2)$ by Lemma~\ref{lem:smoothing}. At each centering, \eqref{eq:parameters} gives a jump of at most $Cs^{-2}(1+z^2)$ in $Z$. In particular, an axial translation of order $s^{-1-\sigma}$ changes the graph and its gradient by at most $Cs^{-1-\sigma}\|u\|_{C^2}=O(s^{-2})$, since $\sigma+\gamma=1$. Both errors are supported on $|z|\leq16s^\beta$ and hence bounded by $A=Cs^{-3/2+2\beta}$.

Set $d(s-h)=0$, solve $d'=2d+A$ between centerings, and increase $d$ by $A$ at each centering. Fix $\eta>0$. The $C^3$ bound on $[s-h,s]$ tends to zero uniformly in the flow, so for large $s$ the zeroth-order error coefficient is below both one and the threshold for $\eta$ in Sun--Xue \cite[Prop.~2.7]{SXII}. Since $L1=1$, $(Z-d)_+$ satisfies the fixed-frame inequality and has no upward jumps. Its Gaussian $L^2$ norm is at most $C/s$. The cutoff error vanishes on $|z|<14r^\beta$, and $d(s)\leq CAe^{2h}=O(s^{-3/2+6\beta})=o(s^{-1+2\beta})$ since $\beta<1/8$. Applying their estimate to $(Z-d)_+$ with decay exponent one and radius $7r^\beta$, then adding $d(s)$, gives
\begin{equation}\label{eq:OU-pointwise}
 \|u(s)\|_{C^1(|z|<6s^\beta)}
       \leq C_\eta s^{-1+2\beta+\eta}.
\end{equation}

Suppose the graph or its $C^3$ bound first fails at $s_*\geq J$ and put $r=s_*-3$. Since $2\beta<1$, the interval $[r-2\log(7r^\beta)-2,r]$ begins after $J-\log J$ for large $J$. On its portion before $J$, the assumed $H_G^1$ and $C^3$ bounds give the same Sun--Xue estimate in the fixed frame, with no rotation or centering errors. Any centering at $J$ itself satisfies \eqref{eq:parameters}. We may therefore fix the frame at $r$ and apply \eqref{eq:OU-pointwise} there on $|z|<6R_r$. The intervening dilations change the final time in this frame to $h_*=3+O(s_*^{-2})$ by \eqref{eq:parameters}. For $1/4\leq h\leq h_*$, Lemma~\ref{lem:cylindrical-pseudolocality} with $\rho=3/4$ gives a small $C^3$ graph on $|z|<6\rho e^{h/2}R_{s_*-3}$. For $0\leq h\leq1/2$ this region lies inside the previously controlled graph, since $6\rho e^{1/4}<16$. At $h_*$ the target has $\xi$-radius $16e^{-h_*/2}R_{s_*}$, which is less than a fixed fraction of $6\rho R_{s_*-3}$ since $16e^{-3/2}<6\rho$. Lemma~\ref{lem:propagate} at $H=h_*\in[1,4]$ therefore bounds the $C^3$ norm of the fixed-frame graph $\widetilde u$ by $Cs_*^{-1+2\beta+\eta}$ on a fixed enlargement of $|z|<16R_{s_*}$. For the graph function $u$ after the intervening centering operations, differentiating the implicit graph transformation formula through order three gives, for an integer $m$ depending only on this derivative order,
\[
 \|u\|_{C^3(|z|<16R_{s_*})}
 \leq C\|\widetilde u\|_{C^3}
       +C(1+R_{s_*})^m s_*^{-1-\sigma}.
\]
Here $s_*^{-1-\sigma}$ bounds the total changes in center, scale and axis over this interval, and the third derivatives enter with bounded coefficients. Choose $\eta>0$ and then $\beta>0$ so that $m\beta<\sigma/2$ and $1-2\beta-\eta>\gamma$. The added geometric term is $o(s_*^{-1+2\beta+\eta})$, so the resulting $C^3$ bound is $C s_*^{-1+2\beta+\eta}$ on $|z|<16R_{s_*}$. This strictly improves the bound at the supposed first failure time, a contradiction. For large $J$, the total loss of radius in \eqref{eq:radius-loss} is less than the margin between $15J^\beta$ and $16J^\beta$.
\end{proof}

\begin{lemma}
\label{lem:reference-seed}
For every sufficiently large $J$, the rotational model of Proposition~\ref{prop:one-neck-reference}, centered at $(0,1)$, satisfies the graphical and norm bounds of Lemma~\ref{lem:extend} on the interval preceding $J$. A dilation $1+O(J^{-2})$ centers the graph at $J$, with strict inequalities in \eqref{eq:bootstrap}.
\end{lemma}

\begin{proof}
Lemma~\ref{lem:neckpinch-comparison}(i) gives the $C^2$ expansion of the model of Proposition~\ref{prop:one-neck-reference} on $B_{\tau^d}$ for some $d>0$. Take $\beta<d$, so that the ball contains the graph on $|z|<64\tau^\beta$ for large $\tau$. For $v=\chi_\tau u$, the Gaussian tail of $H_2$ gives
\begin{equation}\label{rs:H1}
 \left\|v-\frac{H_2}{2\sqrt2\,\tau}\right\|_{H_G^2}\leq C\tau^{-19/10}.
\end{equation}
The same expansion bounds $\|u\|_{C^2}$ by $C\tau^{-1+2\beta}$ on $|z|<64\tau^\beta$. On the preceding unit time interval the substitution $z=e^{h/2}\xi$ removes the axial drift. The interior parabolic estimates in \cite{Lieberman1996} on the smaller interval then give
\begin{equation}\label{rs:coarse}
 \|u(\tau)\|_{C^4(|z|<32\tau^\beta)}\leq C\tau^{-1+2\beta}.
\end{equation}

By symmetry, $v=bH_2+a+w$ with $a$ constant in space and $w\in E_+$. The $C^2$ expansion bounds $|v|+|Dv|+|D^2v|$ by $C(1+z^2)/\tau$, so the Gaussian mean of $\mathcal Q(v)$ is $O(\tau^{-2})$. Projecting the cutoff equation onto constants gives $a'-a=O(\tau^{-2})$, including the exponentially small cutoff error. Since $a\to0$ by \eqref{rs:H1}, integration from infinity gives $a=O(\tau^{-2})$.

By \eqref{rs:H1}, $b=(2\sqrt2\,\tau)^{-1}+O(\tau^{-19/10})$ and $\|w\|_{H_G^1}=O(\tau^{-19/10})$. Lemma~\ref{lem:weighted-IFT} and symmetry therefore give a dilation $1+O(J^{-2})$ as the only correction. By the centering estimates in the proof of Lemma~\ref{lem:center}, the inequalities in \eqref{eq:bootstrap} hold strictly after this dilation if we first choose $C_2$ sufficiently large and then take $J$ large. Estimates \eqref{rs:H1} and \eqref{rs:coarse}, together with the $C^2$ expansion on the larger interval, give the bounds before $J$ since $1-2\beta>\gamma$.
\end{proof}

\subsection{Uniform backward isolation}
\begin{proposition}
\label{prop:uniform-isolation}
Fix $\alpha,\Lambda,\tau_b,c,C>0$. Let $\{K_t\}$ be an $\alpha$-noncollapsed level-set flow with entropy at most $\Lambda$ and a singularity $(q,S)$. Choose coordinates centered at $q$ in which the axis of the limiting round cylinder at $(q,S)$ is the $z$-axis. Put $s=S-t$ and $\tau=-\log s$. Suppose, for $\tau\geq\tau_b$, the normal graph of $u$ on $|z|<c\tau^\beta$ represents the entire boundary $s^{-1/2}(\partial K_t-q)$ in $B_{c\tau^\beta}(0)$, with the rescaled domain on the side containing the axis, and satisfies
\begin{equation}\label{ui:H}
 \|u\|_{C^3(|z|<c\tau^\beta)}\leq C\tau^{-\gamma},
 \qquad \left\|u-\frac{\sqrt2}{4\tau}H_2\right\|_{L_G^2(|z|<c\tau^\beta)}\leq C\tau^{-1-\sigma}.
\end{equation}

There is $r_*>0$, depending only on $\alpha,\Lambda,\tau_b,c,C$, such that $(q,S)$ is the only singularity in $B_{r_*}(q)\times[S-r_*^2,S]$.
\end{proposition}

\begin{proof}
We may take $c\leq1$. Compare the rescaled flow with the fixed rotational model of Proposition~\ref{prop:one-neck-reference}, centered at its neckpinch and rescaled at the same time $\tau$. By \eqref{rs:H1} and \eqref{ui:H}, the two graph functions differ by $O(\tau^{-1-\sigma})$ in $L_G^2$ on $|z|<c\tau^\beta$. On $B_{c\tau^\beta/2}$ the model graph function tends to zero in $C^4$ by \eqref{rs:coarse}. For $u$, apply Lemma~\ref{lem:propagate} on the preceding unit interval with $R=e^{-1/2}c(\tau-1)^\beta$ and $\rho=3/4$. The terminal interval in $z$ has radius $3c(\tau-1)^\beta/4>c\tau^\beta/2$, so the $C^3$ bound in \eqref{ui:H} gives the same convergence. On $B_{c\tau^\beta/3}$, interior regularity bounds the derivatives of the model's unit normal through order four. Since both cylinder graph functions tend to zero in $C^4$, normal projection expresses the rescaled boundary on this smaller ball as a normal graph over the model whose height tends to zero in $C^4$, uniformly over the family. The absolute height at the projected point is bounded by the absolute difference of the two cylinder graph functions at the original cylinder coordinate. The projection has bounded Jacobian and changes the Gaussian exponent by $O((1+\tau^\beta)\tau^{-\gamma})$. Thus the height over the model has Gaussian $L^2$ norm $O(\tau^{-1-\sigma})$ on $B_{c\tau^\beta/3}$, which contains $B_{\tau^{\beta/2}}$ for all sufficiently large $\tau$.

Fix a sufficiently small $\delta>0$ and set $K=(c/4)(3T_0)^{\beta-1/2}$. For $T_0\leq\tau\leq3T_0$ we have $K\sqrt\tau\leq c\tau^\beta/4$, while $K\sqrt{T_0}=(c/4)3^{\beta-1/2}T_0^\beta\to\infty$ as $T_0\to\infty$. The uniform $C^4$ convergence on $B_{c\tau^\beta/3}$ and the uniform Gaussian $L^2$ bound $O(\tau^{-5/4})$ on $B_{\tau^{\beta/2}}$ therefore give a single $T_0$ for which the height over the model has $C^4$ norm at most $\delta$ on $B_{K\sqrt\tau}$ throughout $[T_0,3T_0]$ and Gaussian $L^2$ norm less than $\tau^{-9/8}$ on $B_{\tau^{\beta/2}}$ for every $\tau\geq T_0$. Lemma~\ref{lem:neckpinch-comparison}(ii) now gives the $C^4$ bound $\delta$ on $B_{2\sqrt\tau}$ for all $\tau\geq T_1$ with a common $T_1$. Since $\delta$ was arbitrary, this norm tends uniformly to zero.

In the ambient ball $B_{\sqrt\tau}(0)$, the rescaled boundary therefore agrees with a cylinder graph of uniformly bounded radius. In physical coordinates the axial positions $z=\pm\sqrt\tau/2$ lie at signed distances $\pm h(\tau)/2$ from $q$ along the axis, where $h(\tau)=e^{-\tau/2}\sqrt\tau$ decreases to zero for $\tau>1$. By Lemma~\ref{lem:neckpinch-comparison}(i) with $V=\pm1/2$, the model translated to either position converges locally in $C^1$ to a round cylinder of radius $\sqrt{2+1/4}=3/2$. The uniform $C^4$ estimate gives the same limit for our flows. Measuring space in units of $\sqrt s$ and time in units of $s$ from $S-s$, the cylinder becomes extinct at $(3/2)^2/2=9/8$. Since $9/8>1$, the pseudolocality theorem of Sun--Xue \cite[Thm.~2.5(1)]{SXI}, applied after dilation to the standard cylinder, gives smoothness through $S$ on balls of any fixed rescaled radius $R$ about both centers for $\tau\geq\tau_R$, with $\tau_R$ uniform over the family.

Following the covering argument of Sun--Xue \cite[Sec.~7.1]{SXI}, choose $R$ larger than twice the bound for the rescaled cross-sectional radii, and $\tau_a>1$ so that the preceding regularity holds on balls of rescaled radius $R$ for $\tau\geq\tau_a$. Take $0<r_*<\min\{h(\tau_a)/2,e^{-\tau_a/2}\}$. For $0<|x-q|<r_*$, choose $\tau>\tau_a$ with $h(\tau)=2|x-q|$ and put $s=e^{-\tau}$. For $S-r_*^2\leq t\leq S-s$, monotonicity of $h$ puts $x$ inside the ball where $\partial K_t$ is a smooth graph. Write $z$ for the axial coordinate of $x-q$ and $\rho$ for its distance from the axis. If $x\in\partial K_t$ for $S-s\leq t\leq S$, then $K_t\subset K_{S-s}$ gives $x\in K_{S-s}$, so $\rho\leq R\sqrt s/2$. The two ball centers have axial coordinates $\pm|x-q|$ by the choice of $\tau$. The axial distance to the nearer center is $|x-q|-|z|\leq\rho$, so the distance to that center is at most $\sqrt2\,\rho<R\sqrt s$. The central cross-section is smooth before $S$, which leaves only $(q,S)$.
\end{proof}

\section{Selection of the parameters}\label{sec:selection}
\begin{proposition}
\label{prop:local-package}
Suppose that, in coordinates rescaled about $(0,1)$, a flow satisfies the hypotheses of Lemma~\ref{lem:extend} at $j=J$ and the rescaled domain lies on the side of the graph containing the axis. For every integer $J$ sufficiently large in terms of the constants in these hypotheses, the flow has a nondegenerate neckpinch at $(q,S)$ satisfying
\[
 |q|\leq C e^{-J/2}J^{-1-\sigma},\qquad |S-1|\leq C e^{-J}J^{-2},
\]
where $C$ is independent of $J$ and of the flow. For each fixed $J$, there is a common $r_*>0$ such that $(q,S)$ is the only singularity in $B_{r_*}(q)\times[S-r_*^2,S]$.
\end{proposition}
\begin{proof}
Lemmas~\ref{lem:center} and~\ref{lem:extend} propagate \eqref{eq:bootstrap} at each integer $j\geq J$ and give the graphical $C^3$ bound on $|z|<cj^\beta$ throughout each step. Taking reciprocals in \eqref{eq:discrete-b} and summing gives
\[
 b_{j+1}^{-1}-b_j^{-1}=\kappa+O(j^{-\sigma}),\qquad b_j=\frac1{\kappa j}+O(j^{-1-\sigma}).
\]
For fixed $J$, the bound on $b_J^{-1}$ in \eqref{eq:bootstrap} makes the error in $b_j$ uniform over the family.

Let $q_j$ be the spatial center at the start of step $j$ and $\lambda_{j+1}$ the dilation at its end. In the notation of Section~\ref{sec:local-estimates}, the physical times and scales satisfy
\[
 t_{j+1}=t_j+(1-e^{-1})s_j,\qquad
 s_{j+1}=e^{-1}s_j\lambda_{j+1}^{-2}.
\]
Since $\sum|\lambda_j-1|<\infty$, one has $s_j\asymp s_J e^{-(j-J)}$.  Put $S_j=t_j+s_j$.  Estimates \eqref{eq:parameters} imply
\[
 |S_{j+1}-S_j|\leq Cs_jj^{-2},\qquad
 |q_{j+1}-q_j|\leq C\sqrt{s_j}\,j^{-1-\sigma}.
\]
Consequently $S_j\to S$, $q_j\to q$, and
\begin{equation}\label{eq:physical-tail}
 |S-S_j|\leq Cs_jj^{-2},\qquad
 |q-q_j|\leq C\sqrt{s_j}\,j^{-1-\sigma}.
\end{equation}
The initial slice has time $1-e^{-J}$ and time scale $e^{-J}$ before centering. By the initial centering bound \eqref{eq:parameters}, $|q_J|\leq C e^{-J/2}J^{-1-\sigma}$ and $|S_J-1|\leq C e^{-J}J^{-2}$. Together with \eqref{eq:physical-tail} this gives the asserted bounds on $q$ and $S$. The axis converges with error $O(j^{-1-\sigma})$ by \eqref{eq:parameters} and \eqref{eq:omega}. In the limiting frame, $\tau=-\log(S-t)$ satisfies $\tau=j+O(1)$ on the $j$-th interval, uniformly over these flows. By \eqref{eq:physical-tail}, the change to the limiting frame has relative scale error $O(j^{-2})$ and translation $O(j^{-1-\sigma})$ in rescaled units. Along the rotation to the limiting axis, the cylinder projection stays $C^1$-close to the identity on $|z|<cj^\beta$, with displacement $O(j^{-1-\sigma+\beta})$ and Gaussian weight ratio $\exp(O(j^{-1-\sigma+2\beta}))=1+o(1)$. On a smaller axial interval the surface thus remains a graph throughout the rotation. The $C^3$ bound controls the rotation generators and their first derivatives by $C(1+|z|)$ there. Integrating along the rotation and including the translation and dilation gives a total change of $O(j^{-1-\sigma})$ in $H_G^1$. Throughout the step, the proof of Lemma~\ref{lem:center} gives $b(j+h)=b_j+O(j^{-2})$, $\|a(j+h)\|_{H_G^1}=O(j^{-2})$ and, by \eqref{eq:stable-improve}, $\|w(j+h)\|_{H_G^1}=O(j^{-3/2})$. Thus the flow centered at $(q,S)$ is a normal graph on $|z|<c\tau^\beta$ and has the uniform normal form
\begin{equation}\label{eq:H-normal}
 \left\|u(\tau)-\frac{\sqrt2}{4\tau}H_2\right\|_{H_G^1(|z|<c\tau^\beta)}
       \leq C\tau^{-1-\sigma}.
\end{equation}
The $C^3$ graph bound, the coordinate changes in Lemma~\ref{lem:weighted-IFT}, and \eqref{eq:physical-tail} give $\|u\|_{C^3}\leq Cj^{-3/4}+Cj^{-1-\sigma+\beta}\leq Cj^{-3/4}$ on the smaller interval $|z|<c\tau^\beta$ in the limiting frame. In the corresponding ball the graph is the entire boundary, and the domain contains the axis.

Points on the central cross-sections at times $t_j$ converge to $(q,S)$ and have curvature at least $c s_j^{-1/2}$. Since the spacetime boundary is closed and $s_j\to0$, $(q,S)$ is a singular boundary point. The rescaled graphs converge locally to the round cylinder. The nonzero quadratic term in \eqref{eq:H-normal} gives nondegeneracy by \cite[Def.~1.5]{SXI}. The $H_G^1$ bound in \eqref{eq:H-normal} implies the $L_G^2$ bound in \eqref{ui:H}, so Proposition~\ref{prop:uniform-isolation} gives the common isolation radius for fixed $J$.
\end{proof}

\begin{lemma}
\label{lem:exterior-exclusion}
Suppose $T_j\downarrow0$, $\xi_j\to0$, and $(x_j,t_j)$ are singular points of the flows $\{M_t^{T_j,\xi_j}\}$ with $0\leq t_j\leq T_j$.  Then after passing to a subsequence there is an $i$ such that
\[
 \frac{x_j-p_i}{\sqrt{T_j}}
 \longrightarrow0,
 \qquad
 \frac{t_j}{T_j}\longrightarrow1.
\]
\end{lemma}

\begin{proof}
The initial curvature bound $|A|\leq C/\sqrt{T_j}$ and $(\partial_t-\Delta)|A|^2\leq2|A|^4$ give $t_j/T_j\geq c>0$. Since $K_{t_j}^{T_j,\xi_j}\subset K_0^{T_j,\xi_j}$ and the initial tube has radius at most $C\sqrt{T_j}$, there is $q_j\in\Gamma_{T_j,\xi_j}$ with $|x_j-q_j|\leq C\sqrt{T_j}$.

After passing to a subsequence, either $\min_i|q_j-p_i|/\sqrt{T_j}\to\infty$ or $|q_j-p_i|/\sqrt{T_j}$ stays bounded for some fixed $i$. Lemma~\ref{lem:canonical-pointed-convergence} identifies the flow limits under rescaling about $q_j$ in the first case and $p_i$ in the second as a shrinking cylinder and $M_t^R$, respectively. The rescaled singular points remain in bounded spacetime sets, with times bounded away from zero. By Section~\ref{sec:preliminaries}, every limit of these singular points is singular. The cylinder survives past time $1$ since $a_\infty^2/2>1$ by Proposition~\ref{prop:one-neck-reference}. Hence the limit is the model singularity $(0,1)$.
\end{proof}

\begin{proof}[Proof of Theorem~\ref{thm:main}]
Consider the tori $M_0^{T,\xi}$, $\xi\in Q$, with the common noncollapsing and entropy bounds of Lemma~\ref{lem:global-tube-family}. In the coordinates $(x,t)\mapsto(T^{-1/2}(x-p_i),t/T)$, the $i$-th model pinches at $(\eta\zeta_i,\sigma_i)$, where $\sigma_i=1+\eta\vartheta_i$. For the nonfattening model of Proposition~\ref{prop:one-neck-reference}, Sun--Wang--Xue \cite[Thm.~1.1(i)]{SWXPassing2025} give a compact spacetime neighborhood of its nondegenerate neckpinch $(0,1)$ in $\{t>0\}$ containing no other singularity. Fix $J$ sufficiently large for Lemma~\ref{lem:reference-seed} and Proposition~\ref{prop:local-package}, with the bounds on $q$ and $S-1$ in that proposition small enough to place $(q,S)$ inside this neighborhood.

At the $i$-th neck the change $(x,t)\mapsto((x-\eta\zeta_i)/\sqrt{\sigma_i},t/\sigma_i)$ takes the reference flow to the fixed model. Its rescaled time is $\tau=-\log(1-t/\sigma_i)$. We start the iteration at $\tau=J$, that is, at $t=\sigma_i(1-e^{-J})$ with center $\eta\zeta_i$ and time scale $\sigma_i e^{-J}$ in the coordinates at $p_i$. For fixed $J$, Lemma~\ref{lem:canonical-pointed-convergence} gives smooth convergence to this model on $|z|\leq32\tau^\beta$, $J-\log J\leq\tau\leq J$, together with convergence of the enclosed domains. The convergence is uniform for all $i$ and $\xi\in Q$, since every sequence in $Q$ has a convergent subsequence and there are only finitely many necks. Lemma~\ref{lem:reference-seed} thus gives the bounds before centering in Lemma~\ref{lem:extend} for the whole family, with the domain on the side containing the axis.

At $J$ the two infinitesimal axis rotations have invertible projection onto the rotational modes at the model. The implicit function theorem gives a continuous choice of rotation making $P_{\rm rot}(\chi_Ju)=0$. For this fixed $J$, taking $T$ sufficiently small makes the rotation at most $J^{-2-\sigma}$ and gives the coefficient bounds of Lemma~\ref{lem:weighted-IFT} with the model constants enlarged by a fixed factor. Including this rotation, the initial centering then satisfies \eqref{eq:parameters} uniformly on $Q$. By continuity, the strict inequalities in \eqref{eq:bootstrap} for the centered model in Lemma~\ref{lem:reference-seed} persist for every $i$ and $\xi$. Proposition~\ref{prop:local-package} now applies throughout $Q$. For this fixed $J$, its backward isolation radius becomes $\sqrt{\sigma_i}\,r_*$ in the $T^{-1/2}$-rescaled coordinates. Since $\sigma_i>1/2$, the radius $r_*/\sqrt2$ works for all $i$, $\xi\in Q$ and sufficiently small $T$.

Define $G_T:Q\to(\R^3\times\R)^n$ by taking the $i$-th component of $G_T(\xi)$ to be $(q,S-1)$, where $(q,S)$ is the neckpinch obtained from $M_0^{T,\xi}$ in the rescaled coordinates at $p_i$. For fixed $T$, the initial domains depend smoothly on $\xi$, so level-set stability gives convergence of the flows as $\xi'\to\xi$. The convergence is smooth on the regular graphical regions. By the proof of Lemma~\ref{lem:center}, the derivative of the two rotational coefficients with respect to tilting the axis is $I+O(j^{-1})$ throughout step $j$. The implicit function theorem then gives local choices of axis that agree by uniqueness, so the axis chosen at the start of the step extends continuously in time and in $\xi$ up to the endpoint. Applying Lemma~\ref{lem:weighted-IFT} at the endpoint gives, by induction, continuous dependence on $\xi$ of the center and time estimate at each finite step. The time scales tend to zero uniformly on $Q$, so \eqref{eq:physical-tail} gives uniform convergence of the centers and time estimates to the singular points and times. Hence $G_T$ is continuous.

Let $T_k\downarrow0$ and $\xi_k\to\xi_*$ in $Q$. At each neck, undo its model translation and dilation for $\xi_k$. By Proposition~\ref{prop:local-package} and the choice of $J$, the transformed singular points lie in the fixed model's compact isolating neighborhood. Lemma~\ref{lem:canonical-pointed-convergence} gives convergence of the spacetime boundaries there and smooth convergence near every regular point of the model. Thus every limit of these singular points is the model neckpinch $(0,1)$. Restoring the translations and dilations gives $G_{T_k}(\xi_k)\to\eta\xi_*$. By compactness of $Q$,
\[
 \sup_{\xi\in Q}|G_T(\xi)-\eta\xi|_\infty\longrightarrow0
 \qquad\text{as }T\downarrow0.
\]
For small $T$ this error is less than $\eta$, so the straight-line homotopy from $G_T$ to $\eta\,\mathrm{Id}$ avoids zero on $\partial Q$. Hence $G_T$ has degree one at zero and has a zero $\xi_T\in\operatorname{int}Q$. Every such zero satisfies $\eta|\xi_T|_\infty\leq\sup_{\xi\in Q}|G_T(\xi)-\eta\xi|_\infty$, so $\xi_T\to0$. The corresponding necks pinch at
\begin{equation}\label{dg:matched}
 (p_i,T),\qquad i=1,\ldots,n.
\end{equation}

Suppose that along a sequence $T_j\downarrow0$ the selected flow has a singular point $(y_j,s_j)$ with $0<s_j\leq T_j$ distinct from those in \eqref{dg:matched}. Since $\xi_{T_j}\to0$, Lemma~\ref{lem:exterior-exclusion} gives a subsequence along which the singular points rescaled about some fixed $p_i$ converge to $(0,1)$. Their rescaled times are at most $1$, so they eventually lie in the common backward cylinder of Proposition~\ref{prop:local-package}, whose only singularity is $(0,1)$. This contradicts the choice of $(y_j,s_j)$. Hence for every sufficiently small $T>0$, the flow first becomes singular at time $T$, precisely at $p_1,\ldots,p_n$.
\end{proof}

\bibliographystyle{acm}
\bibliography{main}
\end{document}